\documentclass[smallextended,envcountsect,]{svjour3}
\smartqed
\usepackage{graphicx}
\usepackage[dvipsnames]{xcolor}

\journalname{Numerical Algorithms}
\providecommand{\doi}[1]{\href{https://doi.org/#1}{\nolinkurl{https://doi.org/#1}}}
\usepackage[colorlinks,linkcolor=blue,citecolor=blue,urlcolor=black]{hyperref}

 \usepackage[numbered]{mcode}
\usepackage{times,amsmath,amsbsy,amssymb,amscd,mathrsfs, bm}
\usepackage{algorithm}
\usepackage{algpseudocode}
\usepackage{subcaption, booktabs, multicol, multirow, makecell}
\newcounter{mnote}
\let\oldmarginpar\marginpar
\renewcommand\marginpar[1]{\-\oldmarginpar[\raggedleft\footnotesize #1]%
	{\raggedright\footnotesize #1}}
\newcommand{\dd}{\,{\rm d}}
\newcommand{\inprd}[1]{\left\langle {#1} \right\rangle}
\newcommand{\mV}{\mathcal V}
\newcommand{\mQ}{\mathcal Q}

\newcommand{\IQ}{\mathcal{I}_{\mathcal Q}}

\newcommand{\IQinv}{\mathcal{I}_{\mathcal Q}^{-1}}

\newcommand{\Imu}{\mathcal P}

\definecolor{algorange}{RGB}{230,159,0}      % orange
\definecolor{algskyblue}{RGB}{86,180,233}    % sky blue
\definecolor{alggreen}{RGB}{0,158,115}       % bluish green
\definecolor{algvermillion}{RGB}{213,94,0}   % vermillion
\definecolor{algpurple}{RGB}{204,121,167}    % reddish purple
\definecolor{algdarkred}{RGB}{119,1,58}    % dark red

\newcommand{\LC}[1]{\textcolor{black}{#1}}

\newcommand{\revision}[1]{\textcolor{black}{#1}}

\DeclareMathOperator{\sym}{sym}
\begin{document}

\title{Accelerated Transformed Primal Dual Methods for \revision{Smooth Convex Optimization Problems with Affine Equality Constraints}}
\titlerunning{Accelerated Transformed Primal Dual Methods}

 % \subtitle{Accelerated Transformed Primal Dual Methods}

\author{Long Chen
		\and
		Jingrong Wei}

\institute{
%Communicated by Yurii Nesterov. \\
%-----------------------------------------\\
Long Chen \at
             University of California, Irvine \\
             Irvine, CA 92697, USA\\
             chenlong@math.uci.edu
           \and
              Jingrong Wei (corresponding author) \at
              The Chinese University of Hong Kong \\
              Shatin, Hong Kong SAR, China\\
             jrwei@math.cuhk.edu.hk \\   
}

\date{}
%The correct dates will be entered by the editor.

\maketitle

\begin{abstract}
Accelerated transformed primal--dual (ATPD) methods are proposed for smooth convex optimization problems with affine equality constraints. The methods combine primal--dual dynamics with a transformation designed to accelerate convergence. In the strongly convex case, ATPD achieves an accelerated linear convergence rate and matches the known lower bounds on first-order oracle complexity in terms of gradient evaluations and matrix--vector products. The analysis distinguishes these two types of computational work when assessing the complexity of the method. For general convex objectives, a homotopy perturbation strategy is introduced, and an accelerated sublinear rate for the modified Lyapunov energy is proved under a uniform boundedness condition on the primal iterates at all homotopy stages. The acceleration mechanism is based on an exponentially stable ATPD flow and a modified discrete Lyapunov function. This continuous-time perspective guides the construction and analysis of the discrete methods. Numerical experiments on policy evaluation and constrained signal reconstruction problems support the theory and illustrate the efficiency of the proposed methods.
\end{abstract}


	%		65L20  	Stability and convergence of numerical methods
		%       90C25  	Convex programming	
		%		65B99   Acceleration of convergence
		%      37N40  	Dynamical systems in optimization and economics
            %   47J25 Iterative procedures involving nonlinear operators
            % 65J15 Numerical solutions to equations with nonlinear operators
            % 47H05 Monotone operators and generalizations

%All acknowledgements should be placed in the back of the paper after Conclusions..

%%%%%%%%%%%%%%%%%%%%%%%%%%%
\section{Introduction}  
%%%%%%%%%%%%%%%%%%%%%%%%%%%

We study convex optimization problems with affine equality constraints
\begin{equation}\label{eq: affine equality constrained optimization system}
\min_{x\in \mathcal V} f(x) \quad \text{subject to} \quad Bx=b,
\end{equation}
where $f$ is \revision{convex} and smooth with a Lipschitz continuous gradient, $B:\mathcal V\to\mathcal Q$ is a surjective linear operator between finite-dimensional Hilbert spaces, and $b \in \mathcal Q$ is given. Problems of the form \eqref{eq: affine equality constrained optimization system} arise in numerical partial differential equations (PDEs)~\cite{Benzi.M;Golub.G;Liesen.J2005}, image processing~\cite{Chambolle_Pock_2016,goldstein2014fast}, distributed optimization~\cite{kia2023distributed,dvinskikh2021decentralized}, and machine learning~\cite{lin2020accelerated,alacaoglu2022convergence}.

We assume first that $f$ is $\mu$-strongly convex and $L$-smooth, with $\mu>0$ and condition number $\kappa(f):=L/\mu$. Under these assumptions, \eqref{eq: affine equality constrained optimization system} has a unique solution $x^\star$. For unconstrained problems $\min f(x)$, gradient descent converges at the linear rate $\mathcal O((1-1/\kappa(f))^k)$, while accelerated gradient methods improve this rate to $\mathcal O((1-1/\sqrt{\kappa(f)})^k)$ after $k$ iterations~\cite{nesterov2013introductory,taylor2023optimal,wright2025optimization}. A main challenge is to extend such acceleration to constrained problems like \eqref{eq: affine equality constrained optimization system} without increasing the per-iteration cost.

By introducing a Lagrange multiplier $\lambda$ for the constraint $Bx=b$, problem \eqref{eq: affine equality constrained optimization system} can be written as the \LC{min-max} problem
\begin{equation}\label{eq: min-max problem}
\min_{x\in\mathcal V}\max_{\lambda \in\mathcal Q} \mathcal L(x,\lambda)
:= f(x) + \langle \lambda, Bx-b\,\rangle,
\end{equation}
whose solution $(x^\star,\lambda^\star)$ is \LC{a saddle point}. The problem is strongly convex in the primal variable but not strongly concave in the dual variable. Thus, acceleration techniques for strongly convex--strongly concave saddle-point problems do not apply directly~\cite{kovalev2022accelerated,jin2022sharper}.

Salim et al.~\cite{salim2022optimal} studied \eqref{eq: affine equality constrained optimization system} from a first-order oracle perspective and established the lower complexity bounds
\begin{equation}\label{eq:optcomplexity}
\nabla f:\  \Omega\left (\sqrt{\kappa(f)}\log\frac{1}{\varepsilon} \right),
\quad
B \text{ or } B^{\intercal}:\ \Omega\left  (\sqrt{\kappa(f)\kappa(BB^{\intercal})}\log\frac{1}{\varepsilon} \right ),
\end{equation}
to reach \revision{a point $x$ such that $\|x-x^\star\|^2\le\varepsilon$}. These bounds give the target complexity for projection-free first-order methods, where the basic operations are gradient evaluations $\nabla f$ and matrix--vector products with $B$ or $B^{\intercal}$.

We propose an {\it accelerated transformed primal--dual} (ATPD) method for~\eqref{eq: min-max problem}. For $\mu>0$, ATPD achieves an accelerated linear rate with contraction factor $(1+c/\sqrt{\kappa(f)})^{-1}$ and matches the optimal first-order oracle complexity~\eqref{eq:optcomplexity} without projections or exact subproblem solvers. Since $\mu$ is often hard to estimate, we also introduce Algorithm~\ref{alg: restart perturbed ATPD} with homotopy perturbation. For general convex objectives, this method achieves an accelerated sublinear rate $\mathcal O(1/N^2)$ for the modified Lyapunov energy after $N$ inner iterations \LC{under a uniform boundedness condition on the primal iterates at all homotopy stages}. 
%Thus the method gives a unified optimal treatment of both the strongly convex and general convex settings.

\LC{We next place our contribution in context}. Uzawa-type methods and their variants~\cite{arrow1958studies,bacuta2009schur} solve \eqref{eq: min-max problem} by alternating primal--dual updates, often with nonlinear subproblem solvers or preconditioning~\cite{cheng2000nonlinear,chen1998global,song2019inexact,Hu.Q;Zou.J2002}. Other primal--dual methods include augmented Lagrangian methods (ALM)~\cite{hestenes1969multiplier,powell1969method}, the alternating direction method of multipliers (ADMM)~\cite{boyd2011distributed}, and primal--dual hybrid gradient methods~\cite{chambolle2011first}.

Several accelerated variants of ALM and ADMM have been proposed~\cite{he2010acceleration,he2022inertial,he2022fast,luo2022primal}. However, these methods require solving subproblems of the form
\[\arg\min_x \bigl\{f(x)+\langle\lambda_k,Bx-b\rangle+\tfrac{\beta}{2}\|Bx-b\|^2\bigr\},\]
for some $\beta>0$. Thus, improved iteration counts do not directly imply better oracle complexity.

Acceleration has also been applied to inexact or linearized ALM and ADMM methods. Linearizing the smooth objective or the augmented term can simplify the subproblems and improve iteration complexity~\cite{xu2017accelerated,ouyang2015accelerated,boct2023fast,huang2013accelerated}. \revision{In recent works \cite{he2026accelerated,he2024accelerated,luo2026accelerated}, accelerated linear rates were obtained for linearized ALM methods.} However, accelerated iteration rates alone do not establish the separate oracle complexity bounds in \eqref{eq:optcomplexity}; the cost of any inner solves must also be accounted for.

Other approaches include projected gradient methods~\cite{BOTSARIS198196,patrascu2018convergence}, inexact projection variants~\cite{guo2025inexactprojectedpreconditionedgradient,devolder2014first,patrascu2018convergence}, and penalty methods~\cite{polyak1966constrained,hestenes1969multiplier,rockafellar1976augmented,srivastava2020nesterov}. These methods either require projections, introduce projection errors, or depend on carefully tuned penalty parameters.

The constraint can also be imposed by adding the indicator function of the affine set $\{x:Bx=b\}$. Accelerated proximal gradient methods can then be applied~\cite{tran2018smooth}, but they require the proximal operator of this indicator function, which is the projection onto the affine set. Therefore, such methods fall outside the projection-free first-order setting considered here.

When $f$ is strongly convex ($\mu>0$), fully split primal--dual methods for~\eqref{eq: min-max problem} can converge linearly without projections~\cite{du2019linear,luke2018globally,chen2023transformed}. Accelerated methods have also been developed for decentralized problems~\cite{dvinskikh2021decentralized} and smooth saddle-point problems with bilinear coupling~\cite{kovalev2022accelerated}. For general affine constraints, Salim et al.~\cite{salim2022optimal} attain the separate optimal oracle complexity bounds in \eqref{eq:optcomplexity} using Chebyshev acceleration. \revision{Related work by Borodich et al.~\cite{borodich2025near} studies near-optimal complexity separation for composite saddle-point problems.} ATPD provides an alternative approach based on a transformed primal--dual flow and a modified Lyapunov function.
%; with an appropriate accelerated inner solver, it attains the optimal complexity bounds in \eqref{eq:optcomplexity} in the strongly convex case. For merely convex objectives ($\mu=0$), the homotopy variant achieves an accelerated $\mathcal O(1/N^2)$ rate for the modified Lyapunov energy after $N$ inner iterations, under the uniform boundedness assumption of Theorem~\ref{thm: global convergence rate}.

The rest of the paper is organized as follows. Section~\ref{sec: Preliminary} reviews background from convex optimization. Section~\ref{sec:TPD} revisits non-accelerated transformed primal--dual methods. Section~\ref{sec: ATPD} presents ATPD and its convergence analysis. Section~\ref{sec:convex} extends the method to general convex problems using homotopy perturbation and discusses an augmented Lagrangian variant. Section~\ref{sec: Numerical Results} presents numerical experiments, and Section~\ref{sec: conclusion} gives concluding remarks.

%%%%%%%%%%%%%%%%%%%%%%%%%%%
\section{Notation and Preliminaries}\label{sec: Preliminary}
%%%%%%%%%%%%%%%%%%%%%%%%%%%

Let $\mV$ be a finite-dimensional Hilbert space with inner product $\inprd{\cdot,\cdot}$ and norm $\|\cdot\|$. Let $\mV^{\prime}$ denote the dual space of continuous linear functionals $\mathcal G:\mV\to\mathbb R$. By the Riesz representation theorem, every element of $\mV^{\prime}$ can be identified with an element of $\mV$. We thus use $\inprd{\cdot,\cdot}$ to denote both the inner product in $\mV$ and the duality pairing between $\mV$ and $\mV^{\prime}$.

For a symmetric positive definite (SPD) linear operator $A:\mV\to\mV$,
$$
\inprd{x,y}_A := \inprd{Ax,y}
$$
defines an inner product with the induced norm $\|x\|_A^2:=\inprd{x,x}_A$. If $A$ is symmetric but not positive definite, we use $\|x\|_A^2$ only as notation for the quadratic form $\inprd{Ax,x}$. With this convention, the identity of squares still holds:
\begin{equation}\label{eq:squares}
2\inprd{x,y}_A
= \|x\|_A^2+\|y\|_A^2-\|x-y\|_A^2
= \|x+y\|_A^2-\|x\|_A^2-\|y\|_A^2.
\end{equation}

Let $\sigma(A)$ be the spectrum of $A$, and let $\rho(A)=\max\{|\lambda|:\lambda\in\sigma(A)\}$ be its spectral radius. If $A$ is SPD, then $\sigma(A)\subset\mathbb R_{>0}$. We define the condition number as
\begin{equation}\label{eq:kappaA}
\kappa(A):= \lambda_{\max}(A)/\lambda_{\min}(A),
\end{equation}
where $\lambda_{\max}(A):=\max_{\lambda\in\sigma(A)}\lambda$ and $\lambda_{\min}(A):=\min_{\lambda\in\sigma(A)}\lambda$. 

For two symmetric operators $D$ and $A$, we write $A < D$ if $D-A$ is SPD and $A \le D$ if $D-A$ is positive semidefinite. The relation $A \le D$ is equivalent to $\lambda_{\max}(D^{-1}A) \le 1$ or $\lambda_{\min}(A^{-1}D) \ge 1$ if the operators are SPD.

For any convex $C^1$ function $f:\mV\to\mathbb R$, the Bregman divergence is
$$
D_f(y,x) := f(y)-f(x)-\langle \nabla f(x),\,y-x\rangle, \quad x,y\in\mV.
$$
Since $f$ is convex, $D_f(y,x)\ge 0$. The Bregman divergence is generally not symmetric. Its symmetrized version is given by
\begin{equation}\label{eq: symmetrized Bregman divergenve}
\langle \nabla f(y)-\nabla f(x),\,y-x\rangle = D_f(x,y) + D_f(y,x).
\end{equation}
Direct calculation yields the following three-term identity~\cite{chen1993convergence}, which generalizes the identity of squares \eqref{eq:squares}:
\begin{equation*}\label{eq: Bregman divergence identity}  
   \langle\nabla f(x)-\nabla f(y), x - z \rangle=D_{f}(x, y) + D_{f}(z, x) - D_{f}(z, y).  
\end{equation*}  

A differentiable function $f$ is $\mu_f$-strongly convex with $\mu_f>0$, denoted by
$f\in\mathcal S_{\mu_f}$, if
\begin{equation}\label{eq:Dfbound}
D_f(x,y) \ge \frac{\mu_f}{2}\|x-y\|^2,
\quad \forall\, x,y\in\mV.
\end{equation}
We write $f\in\mathcal S_{\mu_f,L_f}$ if $f$ is also $L_f$-smooth, i.e., for some $L_f>0$,
\begin{equation}\label{eq:Lfbound}
D_f(x,y) \le \frac{L_f}{2}\|x-y\|^2,
\quad \forall\, x,y\in\mV.
\end{equation}
\revision{Throughout, $\mu_f$ denotes the largest constant for which the lower bound in
\eqref{eq:Dfbound} holds, and $L_f$ denotes the smallest constant for which \eqref{eq:Lfbound} holds. With this convention, the condition number of $f$ is $\kappa(f):= L_f/\mu_f$.
When $f$ is quadratic, this definition agrees with the matrix condition number
in~\eqref{eq:kappaA}.}

We use $\mu_A(f)$, $L_A(f)$, and $\kappa_A(f)$ for the strong convexity constant, Lipschitz constant, and condition number in the metric induced by an SPD operator $A$:
$$
\revision{\mu_{A}(f) \|x-y\|_A^2 \le \inprd{\nabla f(x)-\nabla f(y),\,x-y} \le L_{A}(f)\|x-y\|_A^2, \quad  \forall\, x,y\in\mV.}
$$
\LC{Then $L_f = L_I(f)$, $\mu_f = \mu_I(f)$, and $\kappa(f) = \kappa_I(f)$.}
For two SPD operators $M$ and $A$, let $\mu_{A}(M) := \lambda_{\min} (A^{-1}M)$ and $L_{A}(M) := \lambda_{\max} (A^{-1}M)$.

To show the main idea, we often use the linear model $\nabla f(x)=Ax$, where $A$ is a constant SPD matrix satisfying $\mu_f I \le A \le L_f I$. In the nonlinear case, the Bregman divergence \revision{$D_f(x, y)$} plays the role of the quadratic form \revision{$\frac12\|x-y\|_A^2$}. Its gradient with respect to the first variable,
$$
\revision{\nabla D_f(x, y)=\nabla f(x)-\nabla f(y)},
$$
is the nonlinear version of \revision{$Ax-Ay$}. The inequality, following from \eqref{eq: symmetrized Bregman divergenve} and \eqref{eq:Dfbound},
\begin{equation*}\label{eq:AuDf}
\revision{\inprd{\nabla f(x)-\nabla f(y),\,x-y}
\ge D_f(x,y)+\frac{\mu_f}{2}\|x-y\|^2, \quad  \forall\, x,y\in\mV,}
\end{equation*}
replaces the linear operator inequality $A \ge A/2 + \mu_f I/2$.

\section{Transformed Primal-Dual Methods}\label{sec:TPD}
In this section, we summarize a version of the transformed primal--dual (TPD) flow and algorithms from~\cite{chen2023transformed}. TPD can be viewed as a relaxed inexact Uzawa method.

\subsection{\LC{Transformed Primal--Dual Flow}}
%The saddle point $(x^\star,\lambda^\star)$ of the min--max problem \eqref{eq: min-max problem} satisfies the KKT condition
%\begin{equation}\label{eq:original}
%\begin{pmatrix}
%\nabla f & \quad B^{\intercal} \\
%-B & \quad 0
%\end{pmatrix}
%\begin{pmatrix}
%x^\star \\
%\lambda^\star
%\end{pmatrix}
%=
%\begin{pmatrix}
%0 \\
%b
%\end{pmatrix},
%\end{equation}
%
%Even when $f$ is strongly convex, the operator in \eqref{eq:original} is not strongly monotone because it has a zero block. 
The TPD approach uses a linear transformation to regain strong monotonicity of the saddle operator $(\nabla_x\mathcal L,-\nabla_\lambda\mathcal L)$. Let $\boldsymbol z=(x,\lambda)^{\intercal}$ and define
\begin{align}
\label{eq:A}
\mathcal T(\boldsymbol z)
&:= 
\begin{pmatrix}
I & \quad 0 \\
\tfrac{1}{L_f}B & \quad I
\end{pmatrix}
\begin{pmatrix}
\nabla f & \quad B^{\intercal} \\
-B & \quad 0
\end{pmatrix}
\begin{pmatrix}
x \\
\lambda
\end{pmatrix} \\
&=
\begin{pmatrix}
\nabla f & B^{\intercal} \\
-B+\tfrac{1}{L_f}B\nabla f & \quad \tfrac{1}{L_f}BB^{\intercal}
\end{pmatrix}
\begin{pmatrix}
x \\
\lambda
\end{pmatrix},
\notag
\end{align}
where the nonlinear block $\nabla f$ is applied to the primal argument $x$, and $B^{\intercal}$ is the adjoint of $B$ with respect to $\inprd{\cdot,\cdot}$. \LC{The KKT condition of \eqref{eq: min-max problem} is equivalent to $\mathcal T(\boldsymbol z^{\star}) = (0,-b)^{\intercal}$.} 

The linear case $\nabla f(x)=Ax$ shows the main idea, where $A$ is a constant SPD matrix with $\mu_f I \le A \le L_f I$. Since $\inprd{\mathcal T \boldsymbol z, \boldsymbol z} = \inprd{\sym \mathcal T\, \boldsymbol z, \boldsymbol z}$ for $\sym \mathcal T = \tfrac{1}{2}(\mathcal T +\mathcal T^{\intercal})$, we compute the symmetric part:
\begin{align*}
\sym \mathcal T
&=
%\begin{pmatrix}
%A & \tfrac{1}{2L_f}A B^{\intercal} \\
%\tfrac{1}{2L_f}BA & \quad \tfrac{1}{L_f}BB^{\intercal}
%\end{pmatrix} 
%=
\frac{1}{2L_f}
\begin{pmatrix}
L_f A & AB^{\intercal} \\
BA & \quad BB^{\intercal}
\end{pmatrix}
+
\frac{1}{2}
\begin{pmatrix}
A & 0 \\
0 & \quad \tfrac{1}{L_f}BB^{\intercal}
\end{pmatrix}.
\end{align*}
The first matrix is positive semidefinite as $L_fA \ge A^2$ and
 for any $(x,\lambda)^{\intercal}\in \mV\times\mathcal Q$,
% , we have
$$
\begin{pmatrix}
x  & \lambda
\end{pmatrix}
\begin{pmatrix}
L_fA & AB^{\intercal} \\
BA & \quad BB^{\intercal}
\end{pmatrix}
\begin{pmatrix}
x \\ \lambda
\end{pmatrix}
\geq
\|Ax+B^{\intercal}\lambda\|^2 \ge 0.
$$
The following estimate on the cross term is the nonlinear version.

\begin{lemma}[Lemma 3.1 in~\cite{chen2023transformed}]\label{lem:key}
Suppose $f \in \mathcal{S}_{\mu_{f}, L_{f}}$. 
For any $x_1,x_2\in \mathcal V$ and $\lambda_1, \lambda_2\in \mathcal Q$, we have
\begin{align*}
  & \frac{1}{L_f} \inprd{ \nabla f(x_1)-\nabla f(x_2), B^{\intercal}(\lambda_1 - \lambda_2) }\\
    \geq{} &\frac{\mu_{f}}{2} \|v_1-v_2\|^2 - \frac{1}{2L_{f}}\| \lambda_1 - \lambda_2\|^2_{BB^{\intercal}} - \frac{1}{2} \inprd{  \nabla f(x_1)-\nabla f(x_2), x_1 - x_2 },
\end{align*}
where $v = x + L_f^{-1} B^{\intercal} \lambda$ is a transformed variable.
\end{lemma}

Lemma \ref{lem:key} leads to a strong monotonicity property \LC{for the operator $\mathcal T$}. This property works much as strong convexity does for convex optimization.

\begin{theorem}[Adapted from the strong Lyapunov estimate in~\cite{chen2023transformed}, Section~3]\label{thm: TPD A monotone}
Suppose $f \in \mathcal{S}_{\mu_{f}, L_{f}}$. For any $\boldsymbol z_i = (x_i, \lambda_i)^{\intercal} \in \mV \times \mathcal Q, i=1,2$, the operator $\mathcal T$ defined in \eqref{eq:A} is strongly monotone and satisfies
\begin{equation}\label{eq:monotone}
\begin{aligned}
&\inprd{\mathcal T(\boldsymbol z_1)-\mathcal T(\boldsymbol z_2),\boldsymbol z_1-\boldsymbol z_2}\\
&\qquad\geq\frac12\inprd{\nabla f(x_1)-\nabla f(x_2),x_1-x_2}
+\frac{1}{2L_f}\|\lambda_1-\lambda_2\|_{BB^{\intercal}}^2.
\end{aligned}
% (\mathcal T(x) -  \mathcal T(x^{\star}), x - x^{\star}) \geq \frac{1}{2} \inprd{  \nabla f(u) - \nabla f(u^{\star}), u - u^{\star}} + \frac{1}{2L_f}\| p - p^{\star}\|_{BB^{\intercal}}^2.
\end{equation}
%where $S = L_f^{-1} BB^{\intercal}$ and $\IQ$ is an SPD operator approximating $S$. 
\end{theorem}

As $\mathcal T$ is strongly monotone, we consider a simpler version of the TPD flow
\cite[Section~3.1]{chen2023transformed}
\begin{equation}\label{eq: TPD flow}
\boldsymbol z'=-\mathcal P^{-1}(\mathcal T(\boldsymbol z) + (0,b)^{\intercal}),
\quad
\mathcal P=\operatorname{diag}(L_f I, \,\IQ),
\end{equation}
where $\boldsymbol z' = \boldsymbol z'(t)$ denotes the derivative with respect to $t$ and $\IQ$ is an SPD operator approximating $BB^{\intercal}$. Define the Lyapunov function and the vector field in \eqref{eq: TPD flow} as
$$
E(\boldsymbol z):=\frac12\|\boldsymbol z-\boldsymbol z^\star\|_{\mathcal P}^2, \quad \mathcal G(\boldsymbol z) := -\mathcal P^{-1}(\mathcal T(\boldsymbol z) + (0,b)^{\intercal}).
$$
Using the strong monotonicity \eqref{eq:monotone}, we have
$$
\frac{\dd}{\dd t} E(\boldsymbol z(t))
=
\left\langle \nabla E(\boldsymbol z),\boldsymbol z'\right\rangle
=
-\left\langle
\boldsymbol z-\boldsymbol z^\star,
\mathcal T(\boldsymbol z)-\mathcal T(\boldsymbol z^\star)
\right\rangle \le -c_L E(\boldsymbol z(t)).
$$
Here
\[
c_L=\min\{1/\kappa(f),\lambda_{\min}(\IQinv BB^{\intercal})/L_f\}.
\]
It follows that $E(\boldsymbol z(t))\leq e^{-c_Lt}E(\boldsymbol z(0))$.

% Under the implicit Euler method $\boldsymbol z_{k+1} - \boldsymbol z_{k} = \alpha \mathcal G(\boldsymbol z_{k+1})$, we use the strong Lyapunov property to yield
% $$
% \mathcal E(\boldsymbol z_{k+1})  - \mathcal E(\boldsymbol z_{k})  \leq \alpha \inprd{\nabla \mathcal E(\boldsymbol z_{k+1}), \mathcal G(\boldsymbol z_{k+1})} \leq - \alpha c_L \mathcal E(\boldsymbol z_{k+1}),
% $$
% which gives the linear convergence $\mathcal E(\boldsymbol z_{k+1})\leq (1+\alpha c_L)^{-1}\mathcal E(\boldsymbol z_{k})$. }

%with preconditioner $\IV = L_f I$:
%\begin{equation}\label{eq: TPD flow}
%\begin{aligned}
%x' &= -\tfrac{1}{L_f}\bigl[\nabla f(x)+B^{\intercal}\lambda\bigr],\\
%\lambda' &= -\IQ^{-1}\bigl[\tfrac{1}{L_f}B\bigl(\nabla f(x)+B^{\intercal}\lambda\bigr)-Bx+b\bigr],
%\end{aligned}
%\end{equation}
%where $\IQ^{-1}$ is a linear operator to be defined later. 

\subsection{\revision{TPD Algorithm and Convergence Analysis}}

Discretizing the flow \eqref{eq: TPD flow} in time with step size $\alpha >0$ gives the TPD method \cite[Section 4.5]{chen2023transformed}. Algorithm \ref{alg: TPD} provides an equivalent version for implementation.

%\begin{equation}\label{eq: TPD method}
%\begin{aligned}
% \frac{x_{k+1} - x_k}{\alpha} &= -\tfrac{1}{L_f}\left [\nabla f(x_k) + B^{\intercal} \lambda_k\right ],\\
%\frac{\lambda_{k+1} - \lambda_k}{\alpha} &= - \IQ^{-1}\left [ \tfrac{1}{L_f}B(\nabla f(x_{k+1}) + B^{\intercal} \lambda_k) - B(2x_{k+1} - x_k) + b\right ].
%\end{aligned}
%\end{equation}
%Noted that $Bx$ is approximated by $B(2x_{k+1} - x_k)$ using an over-relaxation technique common in accelerated methods~\cite{chambolle2011first,chambolle2016ergodic,wei2024accelerated,heBalancedAugmentedLagrangian2021}. 

\begin{algorithm}
\caption{Transformed Primal-Dual (TPD) Method \cite[Section 4.5]{chen2023transformed}}
\label{alg: TPD}
\linespread{1.325}\selectfont
\begin{algorithmic}[1]
\State \textbf{Input}: $x_0 \in \mathcal V, \lambda_0\in \mathcal Q, L_f > 0, \alpha \in (0, 1/2)$. \\
$\IQ^{-1}$: an SPD approximation of $(BB^{\intercal})^{-1}$ such that $\lambda_{\max}(\IQ^{-1}(BB^{\intercal})) = L_f$. 
\For{$k=0, 1, 2, \dots$}
\State $\displaystyle x_{k+1} = x_k - \tfrac{\alpha}{L_f} (\nabla f(x_{k}) + B^{\intercal}\lambda_{k})$
\State $\displaystyle y_{k+1} = - \tfrac{1}{L_f}(\nabla f(x_{k+1}) + B^{\intercal}\lambda_{k})$
\State $\displaystyle \lambda_{k+1} = \lambda_k + \alpha \IQ^{-1}\left[ B ( y_{k+1} + 2x_{k+1} - x_k ) - b\right] $
\EndFor
\State \textbf{return} $x_{k+1}, \lambda_{k+1}$
\end{algorithmic}
\end{algorithm}

For any fixed $\alpha > 0$, define a shifted Lyapunov function
\begin{equation*}\label{eq:TPD discrete Lyapunov nonlinear}
\mathcal E(\boldsymbol z)
:= \frac{1}{2}\|\boldsymbol z - \boldsymbol z^{\star}\|^2_{\mathcal P - \alpha \mathcal B^{\rm aug}} - \alpha D_f(x^{\star},x), \quad
\LC{\mathcal B^{\rm aug}} :=
\begin{pmatrix}
0 & \quad B^{\intercal} \\
B & \quad L_f^{-1}BB^{\intercal}
\end{pmatrix}.
\end{equation*}
%If $\lambda_{\max}(\IQ^{-1}BB^{\intercal}) \le L_f$, then $\mathcal E(\boldsymbol z) \ge 0$ for $|\alpha| < \tfrac{1}{2}$, and $\mathcal E(\boldsymbol z) = 0$ holds only at $\boldsymbol z = \boldsymbol z^{\star}$~\cite[Lemma 4.5]{chen2023transformed}. 
The following result, adapted from~\cite[Section~4.5]{chen2023transformed}, gives the convergence rate of Algorithm~\ref{alg: TPD}.

\begin{theorem}[Adapted from~\cite{chen2023transformed}, Section~4.5]\label{thm: convergence of G-S TPD}
Suppose $f \in \mathcal{S}_{\mu_{f}, L_{f} }$. Let $\boldsymbol z_k=(x_k,\lambda_k)^{\intercal}$ be generated by Algorithm~\ref{alg: TPD} from $\boldsymbol z_0=(x_0,\lambda_0)^{\intercal}$ with step size $0< \alpha < \frac{1}{2}$. Then
\begin{equation*}
\mathcal E(\boldsymbol z_{k+1})\leq \frac{1}{1+\mu \alpha/2} \mathcal E(\boldsymbol z_k),
\end{equation*}
where $\mu = \min \{1/\kappa(f), 1/\kappa (\IQ^{-1}BB^{\intercal})\}$.
\end{theorem}

Theorem~\ref{thm: convergence of G-S TPD} gives a nonaccelerated convergence guarantee: even when $\kappa(\IQ^{-1}BB^{\intercal})=1$, the resulting gradient-complexity bound is $\mathcal O(\kappa(f)\log(1/\varepsilon))$ for reducing the shifted Lyapunov energy by a factor $\varepsilon\in(0,1)$, with a fixed step size $\alpha\in(0,1/2)$. 
%This motivates the accelerated TPD method.

%Applying $\IQ^{-1}$ can be viewed as an inner iteration that approximates $(BB^{\intercal})^{-1}$. In contrast, projection methods require the exact inverse $(BB^{\intercal})^{-1}$. The error from using $\IQ^{-1}$ is controlled by $\kappa(\IQ^{-1}BB^{\intercal})$.
%

%, which reduces the dependence on $\kappa(f)$ to $\mathcal O(\sqrt{\kappa(f)})$.

%We show that the TPD method is not accelerated.
%\begin{itemize}
%\item If $\IQ = L_f^{-1}\lambda_{\max}(BB^{\intercal}) I$ is a scaled identity, the complexity of Algorithm~\ref{alg: TPD} is
%$$
%\mathcal O\left ((\kappa(f)+\kappa(BB^{\intercal}))\log\frac{1}{\varepsilon}\right)
%$$
%for gradient evaluations and matrix--vector products. 
%%This matches other non-accelerated split methods~\cite{salim2022dualize}.
%
%\smallskip
%
%\item If $\IQ^{-1}$ is an optimal linear solver for $BB^{\intercal}$ such that $\kappa(\IQ^{-1}BB^{\intercal})=\mathcal O(1)$, with a cost of $\mathcal O(\sqrt{\kappa(BB^{\intercal})})$ matrix--vector products, the overall complexity is
%$$
%\nabla f:\ \mathcal O\left(\kappa(f)\log\frac{1}{\varepsilon}\right),
%\quad
%B \text{ or }B^{\intercal}:\ \mathcal O\left(\kappa(f)\sqrt{\kappa(BB^{\intercal})}\log\frac{1}{\varepsilon}\right).
%$$
%The factor $\mathcal O(\kappa(f))$ is suboptimal compared to \eqref{eq:optcomplexity}. This leads to the accelerated TPD method, which reaches $\mathcal O(\sqrt{\kappa(f)})$ complexity.
%\end{itemize}

%%%%%%%%%%%%%%%%%%%%
\section{Accelerated Transformed Primal-Dual Methods}\label{sec: ATPD}
%%%%%%%%%%%%%%%%%%%%
In this section, we present accelerated transformed primal--dual (ATPD) methods by combining acceleration techniques for unconstrained convex optimization with the transformation idea from the TPD flow. We first assume that $\mu_f > 0$.

\subsection{Flow and Stability}
We introduce an auxiliary variable $y$ by splitting $x$ to accelerate the $\nabla f$ part. The Lagrange multiplier $\lambda$ remains unsplit. This leads to the accelerated transformed primal--dual gradient (ATPD) flow:
\begin{subequations}\label{eq:AG TPD}
\begin{align}
\label{eq:ATPDx}
x' &= y-x,\\
\label{eq:ATPDy}
y' &= \frac{1}{2}(x-y)-\frac{1}{\mu_f}\bigl(\nabla f(x)+B^{\intercal}\lambda \bigr),\\
\label{eq:ATPDp}
\lambda' &= \frac{L_f}{2}\,\IQinv\!\left[By-b-\frac{1}{L_f}B\bigl(\nabla f(x)+B^{\intercal}\lambda\bigr)\right].
\end{align}
\end{subequations}
Let $\mathcal G(x,y,\lambda)$ be the vector field defined by the right-hand side of \eqref{eq:AG TPD}. The stationary point $(x^{\star}, x^{\star}, \lambda^{\star})$ satisfies $\mathcal G(x^{\star}, x^{\star}, \lambda^{\star}) = 0$. 
%We use an auxiliary variable $y$ to speed up convergence; $y$ moves toward $x$ along the flow. 
%While removing $y$ results in a second-order ODE, we can prove exponential stability more easily using the first-order version \eqref{eq:AG TPD}.

We define the Lyapunov function:
\begin{equation}\label{eq: ATPD Lyapunov}
 E(x,y,\lambda)
:= D_f(x,x^\star)
+ \frac{\mu_f}{2}\|y-x^\star\|^2
+ \frac{1}{L_f}\|\lambda-\lambda^\star\|_{\IQ}^2 .
\end{equation}
Because $f$ is $\mu_f$-strongly convex, $ E(x,y,\lambda) \ge 0$, and it equals zero only when $x=y=x^\star$ and $\lambda=\lambda^\star$. Thus, $E(x,y,\lambda)$ is a candidate Lyapunov function.

We first consider the linear case $\nabla f(x)=Ax$.
%, where:
%$$
%\nabla f(x)=Ax,\quad \mu_f I\le A\le L_f I,
%$$
%and $A$ is a constant SPD matrix. 
Then $-\inprd{\nabla  E(x,y,\lambda), \mathcal G(x,y,\lambda)}$ is a quadratic form in $v = (x-x^\star,\,y-x^\star,\,\lambda-\lambda^\star)^{\intercal}$ with matrix:
\begin{align*}
M:=
&\begin{pmatrix}
A & & \\
& \mu_f I & \\
& & \tfrac{2}{L_f}\IQ
\end{pmatrix}
\begin{pmatrix}
I & -I & 0 \\
-\tfrac12 I + \mu_f^{-1}A & \tfrac12 I & \tfrac{1}{\mu_f}B^{\intercal} \\
\tfrac12\IQinv BA & -\tfrac{L_f}{2}\IQinv B & \quad \tfrac12\IQinv BB^{\intercal}
\end{pmatrix} \\
&\quad=
\begin{pmatrix}
A & -A & 0 \\
-\tfrac{\mu_f}{2}I + A & \quad \tfrac{\mu_f}{2}I & B^{\intercal} \\
\tfrac{1}{L_f}BA & -B &\  \tfrac{1}{L_f}BB^{\intercal}
\end{pmatrix}.
\end{align*}
Since $v^{\intercal}Mv = v^{\intercal}\sym(M)v$, where $\sym(M) = (M+M^\intercal)/2$, we compute the symmetric part \LC{which removes the skew-symmetric component}:
\begin{align}\label{eq:symM}
\sym(M)
&=
\begin{pmatrix}
A & -\tfrac{\mu_f}{4}I & \quad \tfrac{1}{2L_f}AB^{\intercal} \\
-\tfrac{\mu_f}{4}I & \tfrac{\mu_f}{2}I & 0 \\
\tfrac{1}{2L_f}BA & 0 & L_f^{-1}BB^{\intercal}
\end{pmatrix} \notag \\
&=
\begin{pmatrix}
\tfrac12 A & \quad -\tfrac{\mu_f}{4}I & 0 \\
-\tfrac{\mu_f}{4}I &\quad \tfrac{\mu_f}{2}I & 0 \\
0 & \quad 0 & \tfrac12 L_f^{-1}BB^{\intercal}
\end{pmatrix}
+\frac{1}{2L_f}
\begin{pmatrix}
L_fA & \quad 0 & \quad AB^{\intercal} \\
0 & \quad 0 & 0 \\
BA & \quad 0 & \quad BB^{\intercal}
\end{pmatrix}.
\end{align}
By Lemma~\ref{lem:key}, the second term in \eqref{eq:symM} is symmetric and positive semidefinite. For the first term in \eqref{eq:symM}, using $A\ge \mu_f I$, we obtain:
\begin{align*}
\begin{pmatrix}
\tfrac12 A & \ -\tfrac{\mu_f}{4}I \\
-\tfrac{\mu_f}{4}I &\  \tfrac{\mu_f}{2}I
\end{pmatrix}
&=
\frac14
\begin{pmatrix}
A & \ -\mu_f I \\
-\mu_f I & \ \mu_f I
\end{pmatrix}
+
\frac14
\begin{pmatrix}
A & 0 \\
0 & \mu_f I
\end{pmatrix} \geq
\frac14
\begin{pmatrix}
A & 0 \\
0 & \mu_f I
\end{pmatrix}.
\end{align*}
For the remaining block in \eqref{eq:symM}, we use $BB^{\intercal}\geq \mu_{\IQ}(B) \IQ$ with $\mu_{\IQ}(B) := \lambda_{\min} (\IQinv BB^{\intercal})$.

As a result, the symmetric part of the quadratic form for $-\inprd{\nabla E, \mathcal G}$ is positive definite. This matrix calculation acts as a guide for the nonlinear analysis below.

\begin{theorem}[Strong Lyapunov property]\label{thm: ATPD flow decay}
Assume $f\in \mathcal{S}_{\mu_{f}, L_{f}}$. Let $\IQ$ be an SPD approximation of $BB^{\intercal}$ and $\mu_{\IQ}(B) := \lambda_{\min} (\IQinv BB^{\intercal})$. For the Lyapunov function~\eqref{eq: ATPD Lyapunov} and the accelerated transformed primal--dual flow vector field $\mathcal G$ defined in~\eqref{eq:AG TPD}, the following strong Lyapunov property holds:
\begin{equation*}\label{eq:strong acc TPD}
\begin{aligned}
- \inprd{\nabla  E(x,y,\lambda), \mathcal G(x,y,\lambda)} 
&\geq c_L  E(x,y,\lambda), \quad c_L = \frac{1}{2} \min \{1, \mu_{\IQ}(B)\}.
\end{aligned}
\end{equation*}
Consequently, for any solution $(x(t),y(t),\lambda(t))$ of \eqref{eq:AG TPD}, we have the exponential stability estimate
$$
E(x(t),y(t),\lambda(t))\leq e^{- c_L t} E(x(0),y(0),\lambda(0)), \quad t>0.
$$
\end{theorem}
%\begin{proof}
%\LC{Similar to the matrix version, by direct calculation, we obtain
%$$
%\begin{aligned}
%- \inprd{\nabla \mathcal E, \mathcal G}  ={}& \inprd{ \nabla f(x) - \nabla f(x^{\star}), x - x^{\star}} - \frac{\mu_f}{2} \inprd{ y - x^{\star}, x-y} \\
%    &+ \tfrac{1}{L_f}\inprd{ \lambda-\lambda^{\star}, B(\nabla f(x) - \nabla f(x^{\star}))} + \tfrac{1}{L_f}\|\lambda - \lambda^{\star}\|_{BB^{\intercal}}^2.
%\end{aligned}
%$$}
\begin{proof}
\LC{Since $\mathcal G(x^\star,x^\star,\lambda^\star)=0$, we replace $\mathcal G(x,y,\lambda)$ by $\mathcal G(x,y,\lambda)-\mathcal G(x^\star,x^\star,\lambda^\star)$, use $Bx^\star=b$ and $\nabla f(x^\star)+B^{\intercal}\lambda^\star=0$, and list the gradient of $E$ and the difference of the ATPD vector field below:
{\small
\begin{center}
\renewcommand{\arraystretch}{1.45}
\begin{tabular}{c c c}
\toprule
 & $\nabla E$ & $\mathcal G(x,y,\lambda)-\mathcal G(x^\star,x^\star,\lambda^\star)$ \\
\midrule
$x$ &
$\nabla f(x)-\nabla f(x^\star)$ &
$y-x^\star-(x-x^\star)$ \\
$y$ &
$\mu_f(y-x^\star)$ &
$\frac12(x-y)-\frac{1}{\mu_f}\bigl(\nabla f(x)-\nabla f(x^\star)+B^{\intercal}(\lambda-\lambda^\star)\bigr)$ \\
$\lambda$ &
$\frac{2}{L_f}\IQ(\lambda-\lambda^\star)$ &
$\frac{L_f}{2}\IQ^{-1}\left[B(y-x^\star)-\frac{1}{L_f}B\bigl(\nabla f(x)-\nabla f(x^\star)+B^{\intercal}(\lambda-\lambda^\star)\bigr)\right]$ \\
\bottomrule
\end{tabular}
\end{center}
}} 
Direct substitution gives
\begin{equation*}\label{eq: strong eq 1}
\begin{aligned}
-\inprd{\nabla E,\mathcal G}
={}&
\inprd{\nabla f(x)-\nabla f(x^\star),x-x^\star}
+\frac{1}{L_f}\|\lambda-\lambda^\star\|_{BB^{\intercal}}^2\\
&+
\frac{1}{L_f}\inprd{\lambda-\lambda^\star,B(\nabla f(x)-\nabla f(x^\star))}-\frac{\mu_f}{2}\inprd{y-x^\star,x-y}.
\end{aligned}
\end{equation*}
By Lemma~\ref{lem:key}, we bound the $(x,\lambda)$ cross term as follows:
\begin{equation*}
\begin{aligned}
&\frac12\inprd{\nabla f(x)-\nabla f(x^\star),x-x^\star}
+\frac{1}{2L_f}\|\lambda-\lambda^\star\|_{BB^{\intercal}}^2 \\
&\qquad
+\frac{1}{L_f}\inprd{\lambda-\lambda^\star,B(\nabla f(x)-\nabla f(x^\star))}
\geq 0.
\end{aligned}
\end{equation*}
%Hence
%\begin{equation*}
%\begin{aligned}
%-\inprd{\nabla\mathcal E,\mathcal G}
%\geq{}&
%\frac12\inprd{\nabla f(x)-\nabla f(x^\star),x-x^\star}
%+ \frac{1}{2L_f}\|\lambda-\lambda^\star\|_{BB^{\intercal}}^2 \\
%&-\frac{\mu_f}{2}\inprd{y-x^\star,x-y} .
%\end{aligned}
%\end{equation*}
For the $(x,y)$ cross term, the identity of squares gives
$$
-\frac{\mu_f}{2}\inprd{y-x^\star,x-y}
=
-\frac{\mu_f}{4}\|x-x^\star\|^2
+\frac{\mu_f}{4}\|y-x^\star\|^2
+\frac{\mu_f}{4}\|x-y\|^2.
$$
Using
$$
\inprd{\nabla f(x)-\nabla f(x^\star),x-x^\star}
=
D_f(x,x^\star)+D_f(x^\star,x)
$$
and
$$
D_f(x^\star,x)\geq \frac{\mu_f}{2}\|x-x^\star\|^2,
$$
we obtain
$$
\frac12\inprd{\nabla f(x)-\nabla f(x^\star),x-x^\star}
-\frac{\mu_f}{4}\|x-x^\star\|^2
\geq
\frac12 D_f(x,x^\star).
$$
Dropping the nonnegative term $\frac{\mu_f}{4}\|x-y\|^2$ and using
$$
BB^{\intercal}\geq \mu_{\IQ}(B)\IQ,
$$
we get
$$
-\inprd{\nabla E,\mathcal G}
\geq
\frac12D_f(x,x^\star)
+\frac{\mu_f}{4}\|y-x^\star\|^2
+\frac{\mu_{\IQ}(B)}{2L_f}\|\lambda-\lambda^\star\|_{\IQ}^2.
$$
Therefore,
$$
-\inprd{\nabla E,\mathcal G}
\geq
c_L E,
\qquad
c_L=\frac12\min\{1,\mu_{\IQ}(B)\}.
$$
Along any solution of the ATPD flow,
$$
\frac{\dd}{\dd t} E(x(t),y(t),\lambda(t))
=
\inprd{\nabla E, (x',y',\lambda')}
=
\inprd{\nabla E,\mathcal G}
\leq
-c_L E(x(t),y(t),\lambda(t)).
$$
Grönwall's inequality yields
$$
 E(x(t),y(t),\lambda(t))
\leq
e^{-c_Lt} E(x(0),y(0),\lambda(0)).
$$
This proves exponential stability. \qed
\end{proof}

\subsection{Algorithm}\label{sec: ATPD method}
Using accelerated overrelaxation (AOR)~\cite{wei2024accelerated} and \LC{implicit--explicit} discretization of the ATPD flow \eqref{eq:AG TPD}, we present an accelerated transformed primal--dual (ATPD) method:
\begin{subequations}\label{eq:ATPD original}
\begin{align}
\label{eq:ATPD1}      \frac{ x_{k+1}-x_k}{\alpha}&= y_k - x_{k+1} ,\\
\label{eq:ATPD2}    \frac{y_{k+1}-y_k}{\alpha} &= \frac{1}{2}(x_{k+1} - y_{k+1}) -  \frac{1}{\mu_f}(2\nabla f( x_{k+1}) - \nabla f( x_{k})+ B^{\intercal}\lambda_k), \\
\label{eq:ATPD3}   \frac{\lambda_{k+1} - \lambda_k}{\alpha}&= \frac{L_f}{2} \IQinv \left [B(2y_{k+1} - y_k)-b - \tfrac{1}{L_f}B  (\nabla f(x_{k+1}) + B^{\intercal} \lambda_{k}) \right ].
\end{align}
\end{subequations}
In \eqref{eq:ATPD2} and \eqref{eq:ATPD3}, we use AOR to approximate $\nabla f(x) \approx 2\nabla f( x_{k+1}) - \nabla f( x_{k})$ and $By \approx B(2y_{k+1} -y_k)$. This choice adds symmetry to help the convergence analysis.

\begin{algorithm}
\caption{Accelerated Transformed Primal-Dual (ATPD) Method}
\label{alg:ATPD}
\linespread{1.325}\selectfont
\begin{algorithmic}[1]
\State \textbf{Input}: $x_0, y_0 \in \mathcal V, \lambda_0\in \mathcal Q, L_f,  \mu_f$. $\IQ^{-1}(r,\delta)$: a fixed linear SPD approximate inverse applied to $r$, satisfying \eqref{eq:rhom} with tolerance $\delta\in(0,1)$.
\State \textbf{Set} $\displaystyle 0<\alpha\leq\min\left\{\sqrt{\frac{1-\delta}{2\kappa(f)}},\frac{1-\delta}{1+\delta}\right\}$
\For{$k=0, 1, 2, \dots$}
\State $\displaystyle x_{k+1} = \tfrac{1}{1+\alpha} (\alpha y_k + x_k)$
\State $\displaystyle y_{k+1/2} = \tfrac{1}{\mu_f}(2\nabla f( x_{k+1}) -\nabla f( x_{k}) + B^{\intercal}\lambda_k)$
\State $\displaystyle y_{k+1} = \tfrac{1}{1+\alpha/2}\left (y_k - \alpha y_{k+1/2} + \frac{\alpha}{2} x_{k+1}\right)$
\State $\displaystyle \tilde x_{k+1} =2y_{k+1} - y_k - \tfrac{1}{L_f} (\nabla f( x_{k+1}) + B^{\intercal} \lambda_{k})$
\State $\displaystyle \lambda_{k+1} =\lambda_k + \tfrac{\alpha L_f}{2(1-\delta)}\IQinv(B\tilde x_{k+1} - b, \delta)$
\EndFor
\State \textbf{return} $x_{k+1}, y_{k+1}, \lambda_{k+1}$
\end{algorithmic}
\end{algorithm}

Algorithm \ref{alg:ATPD} gives an implementation using a scaled inexact inner solver. By storing $\nabla f(x_k)$ and $B^{\intercal} \lambda_k$, each iteration needs one gradient evaluation, two matrix--vector products, and one inner linear solve $\IQ^{-1}$.

\subsection{Convergence Analysis}\label{sec:ATPDconvergence}
%\begin{pmatrix}
%\mu_f I & \quad 0 \\
%0 & \quad \tfrac{2}{L_f}\IQ
%\end{pmatrix}
%, \quad 
%\mathcal B^{\rm sym} = \begin{pmatrix}
%0 & \quad B^{\intercal} \\
%B & \quad 0 
%\end{pmatrix}, \quad 
% \end{pmatrix}.
Let  $\Imu := \operatorname{diag}(\mu_f I, \tfrac{2}{L_f}\IQ)$. 
Define $\boldsymbol z:= (y, \lambda)^{\intercal}$ and $\boldsymbol z^{\star}:= (x^{\star}, \lambda^{\star})^{\intercal}$. Denote the Lyapunov function~\eqref{eq: ATPD Lyapunov} at $(x_k, \boldsymbol z_k)$ by
$$
E_k := E(x_k, \boldsymbol z_k) = D_f(x_{k}, x^{\star})  + \frac{1}{2}\|\boldsymbol z_k-\boldsymbol z^{\star}\|^2_{\Imu}.
$$
We expand the Lyapunov difference using the Bregman divergence:
\begin{equation}\label{eq: E diff}
    \begin{aligned}
 E_{k+1}  - E_{k}  ={}&  \inprd{  \nabla_x  E( x_{k+1}, \boldsymbol  z_{k+1}) , x_{k+1} - x_k } -D_f(x_k, x_{k+1})  \\
& +  \inprd{  \nabla_z  E( x_{k+1}, \boldsymbol  z_{k+1}) , \boldsymbol  z_{k+1} -\boldsymbol z_k}- \frac{1 }{2} \| \boldsymbol  z_{k+1} - \boldsymbol z_k\|_{\Imu}^2. 
\end{aligned} 
\end{equation}

\LC{We write \eqref{eq:ATPD original} as a correction of the implicit Euler scheme
$$
\begin{aligned}
\begin{pmatrix}
x_{k+1}-x_k\\
y_{k+1}-y_k\\
\lambda_{k+1}-\lambda_k
\end{pmatrix}
={}&
\alpha \, \mathcal G(x_{k+1},y_{k+1},\lambda_{k+1})
\\
&+
\alpha
\begin{pmatrix}
0 & -I & 0\\[1mm]
-\nabla f/\mu_f & 0 & B^{\intercal}/\mu_f\\[1mm]
0 & (L_f/2)\IQinv B & \quad \frac12\IQinv BB^{\intercal}
\end{pmatrix}
\begin{pmatrix}
x_{k+1}-x_k\\
y_{k+1}-y_k\\
\lambda_{k+1}-\lambda_k
\end{pmatrix},
\end{aligned}
$$
where the nonlinear operator $\nabla f$ acting on $x_{k+1}-x_k$ denotes $\nabla f(x_{k+1})-\nabla f(x_k)$. Pairing the correction terms with $\nabla E$ produces the symmetric cross terms below.}

Substituting these into \eqref{eq: E diff} gives the identity
\begin{equation}\label{eq:intermediate decay}
\begin{aligned}
 E_{k+1} -  E_k ={}& \alpha \inprd{\nabla  E(x_{k+1}, \boldsymbol z_{k+1}), \mathcal G(x_{k+1}, \boldsymbol z_{k+1})}\\
 & -D_f(x_k, x_{k+1}) - \frac{1 }{2} \| \boldsymbol z_{k+1} - \boldsymbol z_k\|_{\Imu}^2    \\
& - \alpha\inprd{\nabla f( x_{k+1}) -\nabla f(x^{\star}), y_{k+1} -y_k} \\
& - \alpha \inprd{  y_{k+1} -x^{\star},\nabla f( x_{k+1}) -\nabla f(x_k) }\\
&   +   \alpha \inprd{ y_{k+1} - x^{\star}, B^{\intercal}  (\lambda_{k+1} - \lambda_{k}) } + \alpha \inprd{\lambda_{k+1} - \lambda^{\star}, B(y_{k+1} - y_{k})} \\
    &+\frac{1}{L_f}\alpha \inprd{\lambda_{k+1} - \lambda^{\star}, BB^{\intercal}(\lambda_{k+1} - \lambda_k)}.
\end{aligned}
\end{equation}
\revision{We write the linear terms involving the coupling of variables $\boldsymbol z = (y, \lambda)^{\intercal}$ as
 \begin{equation}\label{eq:zcross}
\begin{aligned}
    &   \alpha \inprd{ y_{k+1} - x^{\star}, B^{\intercal}  (\lambda_{k+1} - \lambda_{k}) } + \alpha \inprd{\lambda_{k+1} - \lambda^{\star}, B(y_{k+1} - y_{k})} \\
    &+\frac{1}{L_f}\alpha \inprd{\lambda_{k+1} - \lambda^{\star}, BB^{\intercal}(\lambda_{k+1} - \lambda_k)} \\
   ={}& \alpha(\boldsymbol z_{k+1} - \boldsymbol z^{\star}, \boldsymbol z_{k+1} - \boldsymbol z_k )_{ \mathcal B^{\rm aug}}\\
   ={}&   \frac{\alpha}{2}\| \boldsymbol z_{k+1} -\boldsymbol z^{\star} \|_{  \mathcal B^{\rm aug}}^2 + \frac{\alpha}{2} \| \boldsymbol z_{k+1} -\boldsymbol z_k \|_{  \mathcal B^{\rm aug}}^2 - \frac{\alpha}{2} \| \boldsymbol z_{k} -\boldsymbol z^{\star} \|_{  \mathcal B^{\rm aug}}^2,  
\end{aligned}
\end{equation}
where
\[
\mathcal B^{\rm aug}=
\begin{pmatrix}
0 & B^{\intercal} \\
B & L_f^{-1}BB^{\intercal}
\end{pmatrix}
\]
is symmetric and indefinite.
The quadratic terms in \eqref{eq:zcross} can be merged with $\|\boldsymbol z - \boldsymbol z^\star\|_{\mathcal P}^2$ in the Lyapunov function.
}
%The AOR ensures the last three terms in \eqref{eq:intermediate decay} have a symmetric form that allows expansion via squares:
%\begin{align*}
%	&(\boldsymbol z_{k+1} - \boldsymbol z^{\star}, \boldsymbol z_{k+1} - \boldsymbol z_k )_{ \mathcal B^{\rm aug}} \\
%	&=  \frac{1}{2}\| \boldsymbol z_{k+1} -\boldsymbol z^{\star} \|_{  \mathcal B^{\rm aug}}^2 + \frac{1}{2} \| \boldsymbol z_{k+1} -\boldsymbol z_k \|_{  \mathcal B^{\rm aug}}^2 - \frac{1}{2} \| z_{k} -\boldsymbol z^{\star} \|_{  \mathcal B^{\rm aug}}^2.
%\end{align*}

For the nonlinear cross term of the $(x,y)$ variables, we mimic the expansion of squares in \eqref{eq:zcross} and write
\begin{equation}\label{eq: Bregman cross}
\begin{aligned}
& \inprd{ \nabla f(x_{k+1})-\nabla f(x^{\star}), y_{k+1}-y_k }
+ \inprd{ y_{k+1}-x^{\star}, \nabla f(x_{k+1})-\nabla f(x_k) } \\
={}& \inprd{\nabla f(x_{k+1})-\nabla f(x^{\star}), y_{k+1}-x^{\star}} - \inprd{\nabla f(x_k)-\nabla f(x^{\star}), y_k-x^{\star}} \\
&+ \inprd{\nabla f(x_{k+1})-\nabla f(x_k),  y_{k+1}-y_k}.
\end{aligned}
\end{equation}
\LC{In the linear case $\nabla f(x)=Ax$, the identity \eqref{eq: Bregman cross} reduces to the expansion of the cross term associated with the symmetric matrix
$\begin{pmatrix}
0 & A\\
A & 0
\end{pmatrix}$.}

\LC{The expansion of the cross terms in \eqref{eq:zcross} and \eqref{eq: Bregman cross} can be incorporated into the Lyapunov difference in \eqref{eq:intermediate decay}. Therefore, we define the modified Lyapunov function
\begin{equation}\label{eq: ATPD discrete Lyapunov}
\begin{aligned}
\mathcal E^{\alpha M}(x,\boldsymbol z)
:={}& D_f(x,x^{\star})
+\frac{1}{2}\|\boldsymbol z-\boldsymbol z^{\star}\|_{\Imu-\alpha \mathcal B^{\rm aug}}^2\\
&+\alpha\inprd{\nabla f(x)-\nabla f(x^{\star}),y-x^{\star}},
\end{aligned}
\end{equation}
and
$$
\mathcal E^{\alpha M}_k:=\mathcal E^{\alpha M}(x_k,\boldsymbol z_k).
$$
The discussion above leads to the following identity.
}

\begin{lemma}\label{lem:EalphaM-bridge}
Let $(x_k,\boldsymbol z_k)$ be generated by the ATPD method~\eqref{eq:ATPD original}. Then, for any initial guess and step size $\alpha$, the modified Lyapunov function defined in \eqref{eq: ATPD discrete Lyapunov} satisfies 
\begin{equation}\label{eq: E and EB bridge}
\begin{aligned}
\mathcal E^{\alpha M}_{k+1}-\mathcal E^{\alpha M}_k
={}& \alpha \inprd{\nabla  E(x_{k+1},\boldsymbol z_{k+1}),\mathcal G(x_{k+1},\boldsymbol z_{k+1})}\\
&-\Big[
D_f(x_k,x_{k+1})
+\tfrac{1}{2}\|\boldsymbol z_{k+1}-\boldsymbol z_k\|_{\Imu-\alpha \mathcal B^{\rm aug}}^2 \\
&\qquad
+\alpha\inprd{\nabla f(x_{k+1})-\nabla f(x_k),y_{k+1}-y_k}
\Big].
\end{aligned}
\end{equation}
\end{lemma}

We shall show that, for sufficiently small $\alpha$, $\mathcal E^{\alpha M}(x, \boldsymbol z)$ is a valid Lyapunov function. 
In the linear setting, $\mathcal E_k^{\alpha M}$ is a quadratic form in $(x_k-x^\star,\, y_k-x^\star,\, \lambda_k-\lambda^\star)^{\intercal}$. The matrix associated with $\nabla \mathcal E^{\alpha M}$ admits the decomposition:
\begin{equation*}
    \begin{aligned}
        \begin{pmatrix}
            A & \quad \alpha A & \quad 0 \\
            \alpha A &\quad \mu_f I &\quad -\alpha B^{\intercal} \\
           0 & \quad -\alpha B & \quad \tfrac{1}{L_f} (2\IQ - \alpha BB^{\intercal})
        \end{pmatrix} := M_1 + M_2 + M_3,
    \end{aligned}
\end{equation*}
where 
\begin{equation}\label{eq:EalphaM_mat}
    \begin{aligned}
    M_1 &=  \begin{pmatrix}
            A & \quad \alpha A & \quad 0 \\
            \alpha A &\quad  \theta \mu_f I & \quad 0 \\
           0  & \quad 0 & \quad 0 
        \end{pmatrix}, \quad \theta = 1/(1+L_{\IQ}(B)), \\
      M_2 &=       \begin{pmatrix}
           0 & \quad 0 &  \quad 0 \\
     0    &\quad (1-\theta)\mu_f I  &\quad -\alpha B^{\intercal} \\
         0  & \quad  -\alpha B & \quad \tfrac{1}{L_f}  \IQ
       \end{pmatrix}, \quad M_3 =    \begin{pmatrix}
          0& \quad 0 & \quad 0 \\
          0    &\quad 0 &\quad 0 \\
          0    & \quad 0 & \quad  \tfrac{1}{L_f} (\IQ - \alpha BB^{\intercal})
        \end{pmatrix}.           
\end{aligned} 
\end{equation}
\LC{For sufficiently small $|\alpha|$, the nonzero principal blocks of the matrices $M_1, M_2$ and $M_3$ are SPD by a spectral argument.  When $f$ is not quadratic, $M_1$ becomes nonlinear and the spectral analysis does not apply. We therefore prove the following bound.}
\begin{lemma}\label{lem:cross Breg}
    Suppose $f \in \mathcal S_{\mu_f, L_f}$. For any $x, y, \tilde x, \tilde y \in \mathcal V$ and $|\alpha| \leq \sqrt{\mu_f/L_f}$,
    \begin{equation}\label{eq:uwvz}
        D_f(x,\tilde x) + \frac{\mu_f}{2}\|y-\tilde y \|^2 \pm \alpha  \inprd{\nabla f(x)  -\nabla f(\tilde x), y -\tilde y } \geq 0.
    \end{equation} 
If $|\alpha| <\sqrt{\mu_f/L_f}$, the equality in \eqref{eq:uwvz} holds if and only if $x = \tilde x$ and $y = \tilde y $.    
\end{lemma}
\begin{proof}
By Cauchy--Schwarz, Young's inequality, and the smooth-convex gradient inequality $\|\nabla f(x)-\nabla f(\tilde x)\|^2\leq 2L_fD_f(x,\tilde x)$~\cite{nesterov2013introductory}, we have 
    \begin{equation*}
    \begin{aligned}
           &| \inprd{\nabla f(x) - \nabla f(\tilde x), y -\tilde y}|\\
        \leq~& \frac{1}{2\sqrt{\mu_f L_f}}\| \nabla f(x) - \nabla f(\tilde x)\|^2 + \frac{\sqrt{L_f\mu_f}}{2}\|y-\tilde y\|^2\\
          \leq~& \sqrt{L_f/\mu_f}D_f(x, \tilde x)+ \frac{\sqrt{L_f\mu_f}}{2}\|y-\tilde y \|^2.
    \end{aligned}
    \end{equation*}
Therefore, for $|\alpha| \leq \sqrt{\mu_f/L_f}$,
    \begin{align*}
              &D_f(x, \tilde x) + \frac{\mu_f}{2}\|y-\tilde y \|^2 - \alpha \inprd{\nabla f(x) - \nabla f(\tilde x), y - \tilde y  }\\
            \geq & \left (1- |\alpha| \sqrt{L_f/\mu_f}\right) \left (D_f(x,\tilde x) +\frac{\mu_f}{2}\|y-\tilde y \|^2 \right ) \geq 0. 
    \end{align*}
Since $1- |\alpha| \sqrt{L_f/\mu_f}>0$ for $|\alpha| < \sqrt{\mu_f/L_f}$, the result follows. \qed
\end{proof}

We summarize the discussion in the following lemma.

\begin{lemma}\label{lem:EalphaM-positive}
Let $L_{\IQ}(B):=\lambda_{\max}(\IQ^{-1}BB^{\intercal})$. If
\begin{equation}\label{eq:alpha_bound_EalphaM}
|\alpha| \leq \min \left\{\frac{1}{\sqrt{\kappa(f)(1+L_{\IQ}(B))}}, \frac{1}{L_{\IQ}(B)}\right\},
\end{equation}
then
\begin{equation}\label{eq:Erelation}
    0 \leq \mathcal E^{\alpha M}(x,\boldsymbol z)\leq 2 E(x,\boldsymbol z), \quad \forall~ x \in \mathcal V,\;\boldsymbol z \in \mathcal V \times \mathcal Q.
\end{equation} 
If \eqref{eq:alpha_bound_EalphaM} holds strictly, the first equality in \eqref{eq:Erelation} holds if and only if $x=x^\star$ and $\boldsymbol z=\boldsymbol z^\star$.
Moreover, replacing $\boldsymbol z^\star$ by an arbitrary $\tilde{\boldsymbol z}$ gives, for any $x,\tilde x,y,\tilde y\in\mathcal V$,
$$
D_f(x,\tilde x)
+\frac{1}{2}\|\boldsymbol z-\tilde{\boldsymbol z}\|^2_{\Imu-\alpha \mathcal B^{\rm aug}}
+\alpha\inprd{\nabla f(x)-\nabla f(\tilde x),y-\tilde y}
\geq 0 .
$$
\end{lemma}

\begin{proof}
Set $\theta=1/(1+L_{\IQ}(B))$. The decomposition in \eqref{eq:EalphaM_mat} extends to the nonlinear case:
$$
\begin{aligned}
\mathcal E^{\alpha M}(x,\boldsymbol z)
={}& D_f(x,x^\star)
+\frac{\theta\mu_f}{2}\|y-x^\star\|^2
+\alpha\inprd{\nabla f(x)-\nabla f(x^\star),y-x^\star} \\
&+\frac{1}{2}\|(x,\boldsymbol z)-(x^\star,\boldsymbol z^\star)\|_{M_2+M_3}^2 .
\end{aligned}
$$
By Lemma~\ref{lem:cross Breg}, with $\mu_f$ replaced by $\theta\mu_f$, if $|\alpha|\leq \sqrt{\theta/\kappa(f)}$,
$$
D_f(x,x^\star)
+\frac{\theta\mu_f}{2}\|y-x^\star\|^2
+\alpha\inprd{\nabla f(x)-\nabla f(x^\star),y-x^\star}
\geq 0.
$$
The equality holds only when $x=y=x^\star$ if $|\alpha|< \sqrt{\theta/\kappa(f)}$.

It remains to control the quadratic part. The spectral estimates give that $M_2$ is positive semidefinite under the same bound on $\alpha$, $M_3$ is positive semidefinite if $|\alpha|\leq 1/L_{\IQ}(B)$ and 
$L_f^{-1}(\IQ-\alpha BB^{\intercal})$ is positive definite if the bound is strict. Hence $\mathcal E^{\alpha M}(x,\boldsymbol z)\geq 0$, with equality only at $(x,\boldsymbol z)=(x^\star,\boldsymbol z^\star)$ if \eqref{eq:alpha_bound_EalphaM} holds strictly.

Finally, since the same argument applies to $-\alpha$,
$$
 E(x,\boldsymbol z)
=
\frac{1}{2}\mathcal E^{\alpha M}(x,\boldsymbol z)
+
\frac{1}{2}\mathcal E^{-\alpha M}(x,\boldsymbol z)
\geq
\frac{1}{2}\mathcal E^{\alpha M}(x,\boldsymbol z).
$$
The last bound follows by replacing $(x^\star,\boldsymbol z^\star)$ with $(\tilde x,\tilde{\boldsymbol z})$ in the same argument. \qed
\end{proof}
We now state the linear convergence rate for the ATPD method.

\begin{theorem}[Accelerated linear convergence of the ATPD method]\label{thm: linear convergence of ATPD}
Assume $f\in \mathcal{S}_{\mu_{f}, L_{f}}$ and denote $\mu_{\IQ}(B):=\lambda_{\min}(\IQ^{-1}BB^{\intercal}), L_{\IQ}(B):=\lambda_{\max}(\IQ^{-1}BB^{\intercal})$. Let $(x_k,\boldsymbol z_k)$ be generated by the ATPD method~\eqref{eq:ATPD original} with step size $\alpha$ satisfying
\begin{equation}\label{eq:alpha bound}
0<\alpha\leq\min\left\{\frac{1}{\sqrt{\kappa(f)(1+L_{\IQ}(B))}},\frac{1}{L_{\IQ}(B)}\right\}.
\end{equation}
Then, for the Lyapunov function~\eqref{eq: ATPD discrete Lyapunov} and all $k\geq 0$,
\begin{equation}\label{eq: E decay}
\max\{\alpha c_L E_{k+1},\mathcal E^{\alpha M}_{k}\}
\leq
\left(\frac{1}{1+\tfrac12\alpha c_L}\right)^k
\mathcal E^{\alpha M}_{0},
\end{equation}
where $c_L=\tfrac12\min\{1,\mu_{\IQ}(B)\}$.
\end{theorem}

\begin{proof}
By the identity \eqref{eq: E and EB bridge} and the strong Lyapunov property at $(x_{k+1},\boldsymbol z_{k+1})$,
\begin{equation*}
\begin{aligned}
 \mathcal E^{\alpha M}_{k+1} - \mathcal E^{\alpha M}_k \leq&~ - \alpha c_L  E_{k+1} - \Big [D_f(x_k, x_{k+1}) +\frac{1}{2} \left\|\boldsymbol z_{k+1}-\boldsymbol z_k\right\|^2_{\Imu - \alpha \mathcal B^{\rm aug}}\\
&+\alpha \inprd{ \nabla f( x_{k+1}) -\nabla f(x_k), y_{k+1} - y_k} \Big ].
\end{aligned}
\end{equation*}
The term in brackets is nonnegative by Lemma~\ref{lem:EalphaM-positive}. Hence
$$
\mathcal E^{\alpha M}_{k+1}-\mathcal E^{\alpha M}_{k}
\leq
-\alpha c_L E_{k+1} \leq -\alpha c_L \frac12\mathcal E^{\alpha M}_{k+1}.
$$
Therefore,
$$
\mathcal E^{\alpha M}_{k+1}
\leq
\frac{1}{1+\tfrac12\alpha c_L}\mathcal E^{\alpha M}_{k}.
$$
Iterating gives the decay of $\mathcal E^{\alpha M}_{k}$. Moreover,
$$
\alpha c_L E_{k+1}
\leq
\mathcal E^{\alpha M}_{k}-\mathcal E^{\alpha M}_{k+1}
\leq
\mathcal E^{\alpha M}_{k}.
$$
Combining the two estimates proves \eqref{eq: E decay}. \qed
\end{proof}
\paragraph{\LC{Inexact Inner Solvers.}}
%We relate
%$$
%\mu_{\IQ}(B):=\lambda_{\min}\!\bigl(\IQinv BB^{\intercal}\bigr),
%\qquad
%L_{\IQ}(B):=\lambda_{\max}\!\bigl(\IQinv BB^{\intercal}\bigr)
%$$
%to the accuracy of the inner solver. 
Let $S=BB^{\intercal}$. Given $r$, the inner problem is to approximate the solution of
\begin{equation}
Sq^\star=r.
\end{equation}
Consider a convergent linear iteration
\begin{equation}\label{eq:localsolver}
q_{i+1}=q_i-D^{-1}(Sq_i-r),\qquad q_0=0,\quad i=0,1,\ldots,m-1,
\end{equation}
where $D$ is SPD and $\rho:=\|I-D^{-1}S\|_S<1$. Fix $D$ and the iteration count $m$ independently of the right-hand side and the outer iteration. After $m$ steps, the iteration defines an SPD approximate inverse $\IQinv$ by the relation
$$
I-\IQinv S=(I-D^{-1}S)^m.
$$
Thus $\IQinv S$ is a polynomial in $D^{-1}S$, and
$\IQinv=\sum_{j=0}^{m-1}(I-D^{-1}S)^jD^{-1}$.
Choosing $m$ such that $\rho^m\leq\delta$ ensures the accuracy condition
\begin{equation}\label{eq:rhom}
\|I-\IQinv S\|_S\leq\delta.
\end{equation}
More generally, any fixed linear SPD approximate inverse satisfying
\eqref{eq:rhom} can be used. Since $\IQinv S$ is SPD in the $S$-inner
product, this condition implies, with $S=BB^{\intercal}$,
\begin{equation}\label{eq: IQ accuracy}
1-\delta
\le
\lambda_{\min}(\IQinv BB^{\intercal})
\le
\lambda_{\max}(\IQinv BB^{\intercal})
\le
1+\delta.
\end{equation}
Hence
$$
\kappa_{\IQ}(B)
:=
\frac{\lambda_{\max}(\IQinv BB^{\intercal})}
{\lambda_{\min}(\IQinv BB^{\intercal})}
\leq
\frac{1+\delta}{1-\delta}.
$$

% \JW{\begin{remark}
%     In practice, when solving $Sq^* =r$, we verify a neccessary condition of \eqref{eq:rhom}: the inner iteration stops when the relative residual falls below a set tolerance~$\delta \in (0,1)$, that is, to find $q = \IQinv r$ such that 
% \begin{equation}\label{eq:solver_relative_res_prac}
%     \frac{\|(I - \IQinv S)q\|_{S}}{\|q\|_{S}} = \frac{\|r - Sq \|_{\IQinv S \IQinv}}{\|r\|_{\IQinv S \IQinv}} \leq \delta.
% \end{equation}
% \LC{This is not computable!}
% Furthermore, we can simply use $\ell^2$-norm or infinite norm to measure the relative residual.
% \end{remark}}
\begin{remark}
In practice, the condition \eqref{eq:rhom} is not checked directly, since it involves an operator norm. Instead, the inner iteration is stopped by a computable residual test, for example
\begin{equation*}\label{eq:solver_relative_res_prac}
\frac{\|r-Sq\|}{\|r\|}\leq \delta,
\qquad
\text{or}
\qquad
\frac{(D^{-1}(r-Sq),r-Sq)}{(q,r)}\leq \delta^2,
\end{equation*}
where $\|\cdot\|$ may be the Euclidean norm or the infinity norm, and $D^{-1}$ is the operator in \eqref{eq:localsolver} used to define $\IQinv$. Such tests are practical proxies for \eqref{eq:rhom}. By themselves, they do not imply the spectral bound \eqref{eq: IQ accuracy} without further assumptions. Residual-based stopping generally does not define a fixed linear approximate inverse and is therefore not covered by the present convergence theorem.
\end{remark}
%If $D^{-1}$ is scaled so that $\lambda_{\max}(D^{-1}S)\leq 1$, then the eigenvalues of $\IQinv S$ are
%$$
%1-(1-\lambda)^m,\qquad \lambda\in\sigma(D^{-1}S).
%$$
%Thus $\lambda_{\max}(\IQinv S)\leq 1$ for every $m$, and \eqref{eq:rhom} gives the sharper estimate
%$$
%\kappa_{\IQ}(B)\leq \frac{1}{1-\delta}.
%$$

\paragraph{\LC{Rescaling the Inner Solver.}}
Given a fixed SPD preconditioner $\IQ$, set $\widetilde\IQ^{-1}=c\IQ^{-1}$ with $c>0$. Then $\mu_{\widetilde\IQ}(B)=c\,\lambda_{\min}(\IQ^{-1}BB^{\intercal})$ and $L_{\widetilde\IQ}(B)=c\,\lambda_{\max}(\IQ^{-1}BB^{\intercal})$. The effective contraction factor is proportional to
$$
\min\left\{\frac{1}{\sqrt{\kappa(f)(1+L_{\widetilde\IQ}(B))}},\frac{1}{L_{\widetilde\IQ}(B)}\right\}\min\{1,\mu_{\widetilde\IQ}(B)\}.
$$
It is maximized by scaling so that $\mu_{\widetilde\IQ}(B)=1$. Therefore, using \eqref{eq: IQ accuracy}, line~8 of Algorithm~\ref{alg:ATPD} applies the rescaling
$$
\widetilde\IQ^{-1}=\frac{1}{1-\delta}\IQ^{-1},
$$
which gives $\mu_{\widetilde\IQ}(B)\geq1$ and $L_{\widetilde\IQ}(B)\leq(1+\delta)/(1-\delta)$. In convergence estimates for Algorithms~\ref{alg:ATPD} and~\ref{alg: restart perturbed ATPD}, the Lyapunov functions are formed with the fixed effective preconditioner $\widetilde\IQ=(1-\delta)\IQ$. These estimates lead to the following convergence result.

\begin{corollary}[Accelerated rate with inexact constraint preconditioner]\label{cor:ATPD}
Assume the setting of Theorem~\ref{thm: linear convergence of ATPD}. Suppose that a fixed linear SPD inner solver $\IQ^{-1}$ satisfies \eqref{eq:rhom} for some $\delta\in(0,1)$ with $\frac{1+\delta}{1-\delta} \le \sqrt{\kappa(f)}$. 
Then
$$
\alpha=\sqrt{\frac{1-\delta}{2\,\kappa(f)}}
$$
satisfies \eqref{eq:alpha bound} for $\widetilde\IQ$, and Algorithm~\ref{alg:ATPD} satisfies
$$
\max \left \{\frac{\alpha}{2}  E_{k+1},  \mathcal E^{\alpha M}_{k} \right \}
\leq
\left(\frac{1}{1+\alpha/4}\right)^k
 \mathcal E^{\alpha M}_{0},
\qquad k\geq 0.
$$
\end{corollary}

With a moderate tolerance, such as $\delta=0.5$, we have $\kappa_{\IQ}(B)\leq3$. For $\kappa(f)\geq9$, Corollary~\ref{cor:ATPD} applies directly. For smaller $\kappa(f)$, the step-size bound in Algorithm~\ref{alg:ATPD} gives the same asymptotic dependence. Thus the number of gradient evaluations to reach $\varepsilon$-accuracy is
$$
\nabla f:\quad \mathcal O \left(\sqrt{\kappa(f)}\log\frac{1}{\varepsilon}\right).
$$
We now address the cost of the inner solve to form such $\IQinv$. 

Let $m=\dim(\mQ)$ and $n=\dim(\mV)$. Let ${\rm nnz}(B)$ be the number of nonzero entries of $B$. One matrix--vector product $BB^{\intercal}\lambda$ needs $\mathcal O({\rm nnz}(B))$ operations, while forming $BB^{\intercal}$ needs $\mathcal O \bigl(m\,{\rm nnz}(B)\bigr)$ operations. We consider the following cases:

\begin{enumerate}
\item {\it Classical iterative solvers.}
For a suitably scaled Richardson iteration, $\rho=1-\Theta(1/\kappa(BB^{\intercal}))$, and the inner iteration needs $\mathcal O(\kappa(BB^{\intercal})\log(1/\delta))$ steps. Since each step costs $\mathcal O({\rm nnz}(B))$ operations, the total cost of the inner solves is
$$
\mathcal O \left({\rm nnz}(B)\,\sqrt{\kappa(f)}\,\kappa(BB^{\intercal})\log \frac{1}{\delta}\log\frac{1}{\varepsilon}\right).
$$
Richardson iteration gives a linear dependence on $\kappa(BB^{\intercal})$, rather than the square-root dependence in the optimal matrix--vector complexity bound. An accelerated inner solver or a suitable preconditioner can improve this dependence.

\smallskip
\item {\it Accelerated iterative solvers.}
For Chebyshev iteration started from zero with a fixed degree $m$ and spectral bounds for $S=BB^{\intercal}$, the error polynomial satisfies~\cite[Section~6.3]{salim2022optimal}
$$
\|I-\IQinv S\|_S\leq
2\left(\frac{\sqrt{\kappa(S)}-1}{\sqrt{\kappa(S)}+1}\right)^m.
$$
Thus a fixed degree of order $\sqrt{\kappa(S)}\log(2/\delta)$ yields a linear SPD approximate inverse satisfying \eqref{eq:rhom}. For a fixed moderate tolerance $\delta$, the inner iteration is accelerated. The operation complexity for the inner solve is
$$
\mathcal O\left({\rm nnz}(B)\sqrt{\kappa(BB^{\intercal})}\log \frac{1}{\delta}\right).
$$
The total complexity is
$$
\nabla f:\ \mathcal O\left(\sqrt{\kappa(f)}\log\frac{1}{\varepsilon}\right),
\qquad
B \text{ or } B^{\intercal}:\ \mathcal O\left(\sqrt{\kappa(f)\kappa(BB^{\intercal})}\log \frac{1}{\delta}\log\frac{1}{\varepsilon}\right).
$$
For fixed $\delta$, $\log(1/\delta)$ is an absolute constant independent of $\kappa(f)$, $\kappa(BB^{\intercal})$, and $\varepsilon$; \LC{for example, when $\delta=0.1$, we have $\log(1/\delta)<2.4$}. Hence the total oracle complexity matches the lower bounds in~\cite{salim2022optimal}. \LC{In practice, one may use CG as the inner solver, although it is not covered by the present theory because CG is nonlinear as an operator.}

\smallskip
\item {\it Exact solve.}
If $\IQ^{-1}=(BB^{\intercal})^{-1}$, the iteration complexity and the number of gradient evaluations reduce to
$$
\nabla f, (BB^{\intercal})^{-1}: \quad \mathcal O \left(\sqrt{\kappa(f)}\log\frac{1}{\varepsilon}\right).
$$
This matches the optimal rate for unconstrained accelerated first-order methods since it behaves like accelerated projected gradient descent. A Cholesky factorization of $BB^{\intercal}$ costs $\mathcal O(m^3)$ once, plus $\mathcal O \bigl(m\,{\rm nnz}(B)\bigr)$ to form the matrix. This approach is practical only when $m \ll n$; otherwise, the setup cost is too high.

\smallskip
\item {\it Multigrid solvers.}
For numerical PDEs admitting a suitable multilevel hierarchy, suppose a fixed symmetric multigrid cycle defines an SPD approximate inverse with a mesh-independent contraction and $\mathcal O(\mathrm{nnz}(B))$ work per cycle; see, for example, the multigrid construction in Section~\ref{sec: Numerical Results}. The work, measured in units of a fine-grid matrix--vector product, is then
$$
\mathcal O \left(\sqrt{\kappa(f)}\,\log\frac{1}{\delta}\log\frac{1}{\varepsilon}\right).
$$
The setup cost and the availability of such a hierarchy depend on the specific problem.
\end{enumerate}

%%%%%%%%%%%%%%%%%%%%%%%%%%%%
\section{Generalizations to Convex Settings}\label{sec:convex}
%%%%%%%%%%%%%%%%%%%%%%%%%%%%

In this section, we consider the general convex case. We assume that $f$ is convex, coercive, and $L_f$-smooth. Then problem~\eqref{eq: affine equality constrained optimization system} admits at least one solution, although the solution may not be unique. We fix one solution and denote it by $x^\star$. Let $\lambda^\star$ be a Lagrange multiplier satisfying
$\nabla f(x^\star)+B^{\intercal}\lambda^\star=0$.

\subsection{General Convex Functions and Homotopy Perturbation}

For $\varepsilon>0$, consider the perturbed ATPD flow
\begin{subequations}\label{eq: perturbed ATPD}
\begin{align}
x' &= y-x, \\
y' &= \tfrac12(x-y)-\tfrac{1}{\varepsilon}\bigl(\nabla f(x)+B^{\intercal}\lambda\bigr), \\
\lambda' &= \tfrac{L_f}{2}\,\IQinv\!\left[
By-b-\tfrac{1}{L_f}B\bigl(\nabla f(x)+B^{\intercal}\lambda\bigr)
\right].
\end{align}
\end{subequations}

Compared with the ATPD flow, the parameter $\mu_f$ in \eqref{eq:ATPDy} is replaced by the perturbation parameter $\varepsilon$. Let $\boldsymbol z=(y,\lambda)$ and define
\begin{equation}\label{eq: perturbed ATPD Lyapunov}
E(x,\boldsymbol z;\varepsilon)
:=
D_f(x,x^\star)
+\frac{\varepsilon}{2}\|y-x^\star\|^2
+\frac{1}{L_f}\|\lambda-\lambda^\star\|_{\IQ}^2.
\end{equation}

\begin{lemma}[Strong Lyapunov property for the perturbed ATPD flow]\label{lem: perturbed ATPD flow decay}
Assume that $f$ is convex and $L_f$-smooth. Let $\IQ$ be an SPD approximation of $BB^{\intercal}$ and define $\mu_{\IQ}(B):=\lambda_{\min}(\IQinv BB^{\intercal})$. For the Lyapunov function \eqref{eq: perturbed ATPD Lyapunov} and the perturbed flow vector field $\mathcal G$ in \eqref{eq: perturbed ATPD}, we have
\begin{equation}\label{eq:strong acc perturbed TPD}
-\inprd{\nabla E(x,\boldsymbol z;\varepsilon),\mathcal G(x,\boldsymbol z;\varepsilon)}
\geq
c_LE(x,\boldsymbol z;\varepsilon)
-\frac{\varepsilon}{4}\|x-x^\star\|^2,
\end{equation}
where
$c_L=\frac12\min\{1,\mu_{\IQ}(B)\}$.
\end{lemma}
\begin{proof}
The calculation is the same as in Theorem~\ref{thm: ATPD flow decay}, with $\mu_f$ replaced by $\varepsilon$. The only difference is the treatment of the term
$$
-\frac{\varepsilon}{2}\inprd{y-x^\star, x-y}
=
-\frac{\varepsilon}{4}\|x-x^\star\|^2
+\frac{\varepsilon}{4}\|y-x^\star\|^2
+\frac{\varepsilon}{4}\|x-y\|^2 .
$$
In the strongly convex case, the negative term is canceled by
$D_f(x^\star,x)\geq \frac{\mu_f}{2}\|x-x^\star\|^2$.
Here $f$ is only convex, so we only have $D_f(x^\star,x)\geq0$, and the term
$-\frac{\varepsilon}{4}\|x-x^\star\|^2$ remains. Repeating the rest of the proof of Theorem~\ref{thm: ATPD flow decay} gives
$$
-\inprd{\nabla E(x,\boldsymbol z;\varepsilon),\mathcal G(x,\boldsymbol z;\varepsilon)}
\geq
c_LE(x,\boldsymbol z;\varepsilon)
-\frac{\varepsilon}{4}\|x-x^\star\|^2 .
$$
This proves \eqref{eq:strong acc perturbed TPD}. \qed
\end{proof}

Define the Lyapunov function
\begin{equation}\label{eq: perturbed ATPD discrete Lyapunov}
\begin{aligned}
 \mathcal E^{\alpha M}(x,\boldsymbol z;\varepsilon):={}&D_f(x, x^{\star})
+\frac{1}{2}\|\boldsymbol z-\boldsymbol z^{\star}\|_{\mathcal P_{\varepsilon}-\alpha\mathcal B^{\rm aug}}^2\\
&+\alpha\inprd{\nabla f(x)-\nabla f(x^{\star}),y-x^{\star}},
\end{aligned}
\end{equation}
where
$$\mathcal P_{\varepsilon} = 
\begin{pmatrix}
\varepsilon I & \quad 0 \\
0 & \quad \tfrac{2}{L_f}\IQ
\end{pmatrix}, \quad \mathcal  B^{\rm aug} = \begin{pmatrix}
0 & \quad B^{\intercal} \\
B & \quad \tfrac{1}{L_f} BB^{\intercal}
 \end{pmatrix}.$$
For fixed $\varepsilon$, we write $\mathcal E^{\alpha M}_k:=\mathcal E^{\alpha M}(x_k,\boldsymbol z_k;\varepsilon)$.

\begin{lemma}\label{lem:positive_Ealpha_eps}
For convex, $L_f$-smooth $f$, with $L_{\IQ}(B)$ as in Lemma~\ref{lem:EalphaM-positive}, set
$$
\alpha_{\max}:=\min\left\{
\sqrt{\frac{\varepsilon}{L_f(1+L_{\IQ}(B))}},
\frac{1}{L_{\IQ}(B)}
\right\}.
$$
For any $\varepsilon>0$,
\begin{equation}\label{eq:positiveEalphaM}
0\leq \mathcal E^{\alpha M}(x,\boldsymbol z;\varepsilon)\leq 2E(x,\boldsymbol z;\varepsilon),
\quad
\text{if }
|\alpha|\leq \alpha_{\max}.
\end{equation}
Moreover, replacing $\boldsymbol z^\star$ by an arbitrary $\tilde{\boldsymbol z}$ gives, for any $x,\tilde x,y,\tilde y\in\mathcal V$,
$$
D_f(x,\tilde x)
+\frac{1}{2}\|\boldsymbol z-\tilde{\boldsymbol z}\|^2_{\mathcal P_{\varepsilon}-\alpha\mathcal B^{\rm aug}}
+\alpha\inprd{\nabla f(x)-\nabla f(\tilde x),y-\tilde y}
\geq 0.
$$
Furthermore, for any $\theta\in(0,1)$,
\revision{\begin{equation}\label{eq:lowerboundEalphaM}
\mathcal E^{\alpha M}(x,\boldsymbol z;\varepsilon)
\geq
(1-\theta)E(x,\boldsymbol z;\varepsilon),
\quad
\text{if }
|\alpha|\leq
\theta\alpha_{\max}.
\end{equation}}
\end{lemma}

\begin{proof}
The bound \eqref{eq:positiveEalphaM} follows by repeating the proof of Lemma~\ref{lem:EalphaM-positive} with the coefficient of the $y$-quadratic term replaced by $\varepsilon$. The nonnegativity estimate in Lemma~\ref{lem:cross Breg} uses only convexity and $L_f$-smoothness of $f$, and remains valid with any positive coefficient in place of $\mu_f$. The same argument, with $(x^\star,\boldsymbol z^\star)$ replaced by $(\tilde x,\tilde{\boldsymbol z})$, gives the second claim.

It remains to prove \eqref{eq:lowerboundEalphaM}. If $|\alpha|\leq \theta\alpha_{\max}$, then $|\alpha|/\theta\leq \alpha_{\max}$. Since $\mathcal E^{\alpha M}$ is affine in $\alpha$ and $\mathcal E^{0M}=E$, we have
$$
\mathcal E^{\alpha M}(x,\boldsymbol z;\varepsilon)
=
(1-\theta)E(x,\boldsymbol z;\varepsilon)
+
\theta\mathcal E^{\frac{\alpha}{\theta}M}(x,\boldsymbol z;\varepsilon).
$$
By \eqref{eq:positiveEalphaM}, $\mathcal E^{\frac{\alpha}{\theta}M}(x,\boldsymbol z;\varepsilon)\geq 0$. Hence
$$
\mathcal E^{\alpha M}(x,\boldsymbol z;\varepsilon)
\geq
(1-\theta)E(x,\boldsymbol z;\varepsilon).
$$
This proves \eqref{eq:lowerboundEalphaM}.
\end{proof}

\LC{
To establish convergence of the perturbed ATPD method, we first prove a decay estimate under an a priori boundedness assumption on the primal iterates. For a fixed perturbation parameter, we verify this boundedness assumption using an exterior quadratic growth condition.
}

\begin{lemma}[\revision{Perturbed decay under boundedness}]\label{lem:perturbed_decay_bounded}
Assume that $f$ is convex and $L_f$-smooth. Let $\IQ^{-1}$ be an SPD approximation of $(BB^{\intercal})^{-1}$ with
$$
\mu_{\IQ}(B):=\lambda_{\min}(\IQ^{-1}BB^{\intercal})\geq 1,
\qquad
L_{\IQ}(B):=\lambda_{\max}(\IQ^{-1}BB^{\intercal}).
$$
Let $(x_k,\boldsymbol z_k):=(x_k,y_k,\lambda_k)$ be generated by \eqref{eq:ATPD original} with $\mu_f=\varepsilon>0$ and step size
\begin{equation*}\label{eq:alpha_bound_perturbed}
0<\alpha\leq
\min\left\{
\sqrt{\frac{\varepsilon}{L_f(1+L_{\IQ}(B))}},
\frac{1}{L_{\IQ}(B)}
\right\}.
\end{equation*}
Assume that the primal iterates are uniformly bounded, namely
\begin{equation}\label{eq:boundedness_assumption}
\|x_k-x^\star\|^2\leq R^2,
\qquad k\geq 0.
\end{equation}
Then
\begin{equation}\label{eq:perturbed_linear_decay}
\mathcal E^{\alpha M}_{k}
\leq
\left(\frac{1}{1+\alpha/4}\right)^k\mathcal E^{\alpha M}_0
+R^2\varepsilon,
\qquad k\geq 0,
\end{equation}
and
\begin{equation}\label{eq:perturbed_Elinear_decay}
E_{k+1}
\leq
\frac{2}{\alpha}
\left(\frac{1}{1+\alpha/4}\right)^k\mathcal E^{\alpha M}_0
+
\left(\frac{2}{\alpha}+\frac12\right)R^2\varepsilon,
\qquad k\geq 0.
\end{equation}
\end{lemma}

\begin{proof}
Following Section~\ref{sec:ATPDconvergence}, the identity \eqref{eq: E and EB bridge} holds with $\mu_f=\varepsilon$:
\begin{equation}\label{eq: E and EB bridge_eps}
\begin{aligned}
\mathcal E^{\alpha M}_{k+1}-\mathcal E^{\alpha M}_k
={}&
\alpha\inprd{\nabla E(x_{k+1},\boldsymbol z_{k+1}),\mathcal G(x_{k+1},\boldsymbol z_{k+1})} \\
&-
\Bigg[
D_f(x_k,x_{k+1})
+\frac12\|\boldsymbol z_{k+1}-\boldsymbol z_k\|^2_{\mathcal P_\varepsilon-\alpha\mathcal B^{\rm aug}} \\
&\qquad
+\alpha\inprd{\nabla f(x_{k+1})-\nabla f(x_k),y_{k+1}-y_k}
\Bigg].
\end{aligned}
\end{equation}
By Lemma~\ref{lem:positive_Ealpha_eps}, the bracketed term is nonnegative. Hence, using Lemma~\ref{lem: perturbed ATPD flow decay},
\begin{equation}\label{eq:discrete-perturbed-decay}
\begin{aligned}
\mathcal E^{\alpha M}_{k+1}-\mathcal E^{\alpha M}_{k}
&\leq
-\frac{\alpha}{2}E_{k+1}
+\frac{\alpha\varepsilon}{4}\|x_{k+1}-x^\star\|^2 \\
&\leq
-\frac{\alpha}{4}\mathcal E^{\alpha M}_{k+1}
+\frac{\alpha\varepsilon}{4}\|x_{k+1}-x^\star\|^2,
\end{aligned}
\end{equation}
where we used $\mathcal E^{\alpha M}_{k+1}\leq 2E_{k+1}$ from \eqref{eq:positiveEalphaM} with the required bound for the step size. By the boundedness assumption \eqref{eq:boundedness_assumption},
$$
\mathcal E^{\alpha M}_{k+1}
\leq
\frac{1}{1+\alpha/4}
\left(
\mathcal E^{\alpha M}_{k}
+\frac{\alpha\varepsilon}{4}R^2
\right).
$$
Iterating this recursion gives
$$
\mathcal E^{\alpha M}_{k}
\leq
\left(\frac{1}{1+\alpha/4}\right)^k\mathcal E^{\alpha M}_0
+
\frac{\alpha\varepsilon R^2}{4+\alpha}
\sum_{i=0}^{k-1}
\left(\frac{1}{1+\alpha/4}\right)^i.
$$
The geometric sum is bounded by $(4+\alpha)/\alpha$, which proves \eqref{eq:perturbed_linear_decay}.

Finally, \eqref{eq:discrete-perturbed-decay} gives
$$
\frac{\alpha}{2}E_{k+1}
\leq
\mathcal E^{\alpha M}_{k}
-\mathcal E^{\alpha M}_{k+1}
+
\frac{\alpha\varepsilon}{4}\|x_{k+1}-x^\star\|^2
\leq
\mathcal E^{\alpha M}_{k}
+
\frac{\alpha\varepsilon}{4}R^2.
$$
Using \eqref{eq:perturbed_linear_decay} yields \eqref{eq:perturbed_Elinear_decay}.
\end{proof}

\LC{
For a fixed perturbation parameter, we next verify the boundedness assumption \eqref{eq:boundedness_assumption} under an exterior quadratic growth condition.
}

\begin{proposition}[\revision{Boundedness from exterior quadratic growth}]\label{prop:boundedness_exterior_qg}
Assume that $f$ is convex, $L_f$-smooth, and satisfies the exterior quadratic growth condition
\begin{equation}\label{eq:exterior-qg}
D_f(x,x^\star)
\geq
\gamma\|x-x^\star\|^2,
\qquad
\|x-x^\star\|\geq R_0,
\end{equation}
for some $\gamma>0$ and $R_0>0$. Assume that the fixed SPD preconditioner satisfies $\mu_{\IQ}(B)\geq1$. Let $(x_k,\boldsymbol z_k)$ be generated by \eqref{eq:ATPD original} with $\mu_f=\varepsilon$ and step size
\begin{equation}\label{eq:alphatheta}
0<\alpha\leq
\theta
\min\left\{
\sqrt{\frac{\varepsilon}{L_f(1+L_{\IQ}(B))}},
\frac{1}{L_{\IQ}(B)}
\right\},
\end{equation}
for some $\theta\in(0,1)$. If $\varepsilon\leq \gamma$, then the primal iterates are uniformly bounded:
$$
\|x_k-x^\star\|\leq R,
\qquad k\geq 0,
$$
where one may take
$$
R^2
:=
\max\left\{
R_0^2,\,
\frac{1}{\gamma(1-\theta)}
\max\{\mathcal E^{\alpha M}_0,\varepsilon R_0^2\}
\right\}.
$$
\end{proposition}

\begin{proof}
From \eqref{eq:discrete-perturbed-decay}, we have
\begin{equation}\label{eq:discrete-perturbed-decay-for-bound}
\mathcal E^{\alpha M}_{k+1}-\mathcal E^{\alpha M}_{k}
\leq
-\frac{\alpha}{2}E_{k+1}
+\frac{\alpha\varepsilon}{4}\|x_{k+1}-x^\star\|^2.
\end{equation}
We first prove that $\mathcal E^{\alpha M}_k$ is uniformly bounded by induction. Assume that
$$
\mathcal E^{\alpha M}_{k}
\leq
C_E:=
\max\{\mathcal E^{\alpha M}_0,\varepsilon R_0^2\}.
$$
If $\|x_{k+1}-x^\star\|\geq R_0$, then \eqref{eq:exterior-qg} gives
$$
\frac{\alpha\varepsilon}{4}\|x_{k+1}-x^\star\|^2
\leq
\frac{\alpha\varepsilon}{4\gamma}D_f(x_{k+1},x^\star)
\leq
\frac{\alpha\varepsilon}{4\gamma}E_{k+1}
\leq
\frac{\alpha}{4}E_{k+1}.
$$
Therefore, \eqref{eq:discrete-perturbed-decay-for-bound} implies
$$
\mathcal E^{\alpha M}_{k+1}
\leq
\mathcal E^{\alpha M}_{k}
\leq C_E.
$$
If $\|x_{k+1}-x^\star\|<R_0$, then
$$
\frac{\alpha\varepsilon}{4}\|x_{k+1}-x^\star\|^2
\leq
\frac{\alpha\varepsilon}{4}R_0^2.
$$
Using also $\mathcal E^{\alpha M}_{k+1}\leq 2E_{k+1}$, we obtain
$$
\mathcal E^{\alpha M}_{k+1}
\leq
\frac{1}{1+\alpha/4}
\left(
\mathcal E^{\alpha M}_{k}
+\frac{\alpha\varepsilon}{4}R_0^2
\right)
\leq
C_E.
$$
Since $\mathcal E^{\alpha M}_0\leq C_E$, induction gives
$$
\mathcal E^{\alpha M}_k\leq C_E,
\qquad k\geq 0.
$$

Since $\alpha$ satisfies the stricter step-size condition \eqref{eq:alphatheta} with $\theta$, Lemma~\ref{lem:positive_Ealpha_eps} gives
$$
\mathcal E^{\alpha M}_k
\geq
(1-\theta)E_k.
$$
Hence
$$
D_f(x_k,x^\star)
\leq
E_k
\leq
(1-\theta)^{-1}\mathcal E^{\alpha M}_k
\leq
(1-\theta)^{-1}C_E.
$$
If $\|x_k-x^\star\|\geq R_0$, then \eqref{eq:exterior-qg} gives
$$
\|x_k-x^\star\|^2
\leq
\gamma^{-1}D_f(x_k,x^\star)
\leq
\frac{C_E}{\gamma(1-\theta)}.
$$
Inside the ball $\|x_k-x^\star\|<R_0$, the same bound is trivial after taking the maximum with $R_0^2$. This proves the stated bound. \qed
\end{proof}

\begin{remark}[\revision{On weaker boundedness assumptions}]
The proof of Lemma~\ref{lem:perturbed_decay_bounded} only requires the a priori boundedness condition \eqref{eq:boundedness_assumption}. Proposition~\ref{prop:boundedness_exterior_qg} verifies it for a fixed perturbation parameter. The homotopy theorem below separately assumes a uniform bound over all stages; the fixed-parameter result does not by itself establish that bound. Other assumptions may also be used, provided they imply a uniform bound on the primal iterates. Notice that coercivity alone gives existence but does not by itself control the quadratic perturbation term in the Lyapunov estimate.
\end{remark}

\begin{algorithm}[htbp]
    \caption{ATPD with homotopy perturbation.}
    \label{alg: restart perturbed ATPD}
\linespread{1.325}\selectfont
    \begin{algorithmic}[1]
        \State \textbf{Parameters}: Initial values $x_0, y_0 \in \mV$, $\lambda_0 \in \mathcal Q$, $L_f,R>0$, $\theta\in(0,1)$, and initial $\varepsilon_0>0$, with $R$ and $\varepsilon_0$ satisfying the assumptions of Theorem~\ref{thm: global convergence rate}. \\
 $\IQ^{-1}(r,\delta)$: a fixed linear SPD approximate inverse applied to $r$, satisfying \eqref{eq:rhom} with tolerance $\delta\in(0,1)$.
        \State \textbf{Set} $m_0
=
\left(
\frac{4}{\theta}
\sqrt{\frac{2L_f}{1-\delta}}
+
\sqrt{\varepsilon_0}
\right)
\ln\!\bigl(4(1-\theta)^{-1}(R^2+1)\bigr)
\varepsilon_0^{-1/2}
$.
       \For{$k=0, 1, 2, \dots$}
        
        \State $\varepsilon_{k+1} = \varepsilon_k/2, \quad m_{k+1} = \sqrt{2}\, m_k$
        
\State Apply Algorithm \ref{alg:ATPD} with $\mu_f = \varepsilon_{k+1}$, initial values $(x_k, y_k, \lambda_k)$, and step size $\alpha=\theta \sqrt{(1-\delta)\varepsilon_{k+1}/2L_f}$ for $\lceil m_{k+1}\rceil$ iterations to find $(x_{k+1}, y_{k+1}, \lambda_{k+1})$.

        \EndFor
        \State \textbf{return} $(x_{k+1}, y_{k+1}, \lambda_{k+1})$
    \end{algorithmic}
\end{algorithm}

Since the perturbation does not change the equilibrium $x^\star$, we use a decreasing sequence $\{\varepsilon_k\}$ and a homotopy argument to obtain an accelerated decay estimate for the modified Lyapunov energy in terms of the total number of inner iterations, as stated in the following theorem. Algorithm~\ref{alg: restart perturbed ATPD} summarizes the method.

\begin{theorem}[\revision{Global convergence of the homotopy ATPD method}]\label{thm: global convergence rate}
Assume that $f$ is convex and $L_f$-smooth. Suppose that a fixed linear SPD inner solver $\IQ^{-1}$ satisfies \eqref{eq:rhom} for some $\delta\in(0,1)$. All energies in this theorem use the fixed effective preconditioner $\widetilde\IQ=(1-\delta)\IQ$. Let $\theta\in(0,1)$ and $\varepsilon_0>0$ satisfy
$$
\frac{1+\delta}{1-\delta}\leq \sqrt{\frac{L_f}{\varepsilon_0}}.
$$
Define
$$
\varepsilon_k:=\varepsilon_0 2^{-k},
\qquad
\alpha_k:=\theta\sqrt{\frac{(1-\delta)\varepsilon_k}{2L_f}},
\qquad k\geq 0.
$$
Assume that all primal iterates generated during the inner ATPD iterations at every homotopy stage satisfy
\[
\|x-x^\star\|^2\le R^2,
\]
for some constant $R>0$, and enlarge $R$ if necessary so that
$$
R^2\geq
\frac{\mathcal E^{\alpha_0M}(x_0,y_0,\lambda_0;\varepsilon_0)}{\varepsilon_0}.
$$
Then the sequence generated by Algorithm~\ref{alg: restart perturbed ATPD} satisfies
\begin{equation}\label{eq:Ek-ek}
\mathcal E^{\alpha_kM}(x_k,y_k,\lambda_k;\varepsilon_k)
\leq
(R^2+1)\varepsilon_k,
\qquad k\geq 0.
\end{equation}
Moreover, for the nominal count $N_k:=\sum_{i=0}^{k}m_i$, we have
\begin{equation}\label{eq:rate}
\varepsilon_k^{-1/2}
=
CN_k+\frac{\varepsilon_0^{-1/2}}{\sqrt{2}},
\end{equation}
where
$$
C=
\frac{\sqrt{2}-1}
{\sqrt{2}\left(\frac{4}{\theta}\sqrt{2L_f/(1-\delta)}+\sqrt{\varepsilon_0}\right)
\ln(4(1-\theta)^{-1}(R^2+1))}.
$$
The actual count is $\widehat N_k:=\sum_{i=1}^{k}\lceil m_i\rceil\leq N_k+k$. Consequently, the homotopy ATPD method reaches modified Lyapunov energy $\mathcal O(\varepsilon)$ in
$$
\widehat N_k=\mathcal O\left(\sqrt{\frac{L_f}{\varepsilon}}\right)
$$
inner iterations. By \eqref{eq:lowerboundEalphaM}, for $k\geq1$ the fixed error quantities satisfy
$$
D_f(x_k,x^\star)+\frac{1}{L_f}\|\lambda_k-\lambda^\star\|_{\widetilde\IQ}^2
\leq\frac{R^2+1}{1-\theta}\varepsilon_k
=\mathcal O(L_f/\widehat N_k^2).
$$
\end{theorem}

\begin{proof}
We prove \eqref{eq:Ek-ek} by induction. The case $k=0$ follows from the assumption
$$
\mathcal E^{\alpha_0M}(x_0,y_0,\lambda_0;\varepsilon_0)
\leq
R^2\varepsilon_0
\leq
(R^2+1)\varepsilon_0.
$$
Suppose that \eqref{eq:Ek-ek} holds for some $k$. By the construction of $\IQ$, the rescaling in Algorithm~\ref{alg:ATPD}, and the choice of $\alpha_{k+1}$, the assumptions of Lemma~\ref{lem:perturbed_decay_bounded} hold with $\varepsilon=\varepsilon_{k+1}$, $\alpha=\alpha_{k+1}$, and the same bound $R$. Applying Lemma~\ref{lem:perturbed_decay_bounded} during the $\lceil m_{k+1}\rceil$ inner ATPD steps with fixed parameters $(\varepsilon_{k+1},\alpha_{k+1})$ gives
\begin{equation}\label{eq:decay_k+1stage}
\mathcal E^{\alpha_{k+1}M}_{k+1}
\leq
\left(\frac{1}{1+\alpha_{k+1}/4}\right)^{m_{k+1}}
\mathcal E^{\alpha_{k+1}M}(x_k,y_k,\lambda_k;\varepsilon_{k+1})
+
R^2\varepsilon_{k+1},
\end{equation}
where
$$
\mathcal E^{\alpha_{k+1}M}_{k+1}
:=
\mathcal E^{\alpha_{k+1}M}(x_{k+1},y_{k+1},\lambda_{k+1};\varepsilon_{k+1}).
$$

Since $\varepsilon_{k+1}=\varepsilon_k/2$ and $m_{k+1}=\sqrt{2}m_k$, we have
$$
m_{k+1}\sqrt{\varepsilon_{k+1}}
=
m_0\sqrt{\varepsilon_0}.
$$
By the choice of $m_0$ in Algorithm~\ref{alg: restart perturbed ATPD}, we have
\begin{equation}\label{eq:decay_rate_bound}
\left(\frac{1}{1+\alpha_{k+1}/4}\right)^{m_{k+1}}
\leq
\exp\left(-\frac{\alpha_{k+1}m_{k+1}}{4+\alpha_0}\right)
\leq
\frac{1-\theta}{4(R^2+1)}.
\end{equation}

We next compare the Lyapunov functions before and after changing the homotopy parameter. Since $\varepsilon_{k+1}<\varepsilon_k$, and since the same $\alpha_{k+1}$ is used in the first two terms below, Lemma~\ref{lem:positive_Ealpha_eps} gives
\begin{equation}\label{eq:EalphaM_bridge}
\begin{aligned}
\mathcal E^{\alpha_{k+1}M}(x_k,y_k,\lambda_k;\varepsilon_{k+1})
&\leq
\mathcal E^{\alpha_{k+1}M}(x_k,y_k,\lambda_k;\varepsilon_k) \\
&\leq
2E(x_k,y_k,\lambda_k;\varepsilon_k) \\
&\leq
2(1-\theta)^{-1}
\mathcal E^{\alpha_kM}(x_k,y_k,\lambda_k;\varepsilon_k).
\end{aligned}
\end{equation}
Combining \eqref{eq:decay_k+1stage}, \eqref{eq:decay_rate_bound}, \eqref{eq:EalphaM_bridge}, and the induction hypothesis gives
$$
\begin{aligned}
\mathcal E^{\alpha_{k+1}M}_{k+1}
&\leq
\frac{1-\theta}{4(R^2+1)}
2(1-\theta)^{-1}
\mathcal E^{\alpha_kM}(x_k,y_k,\lambda_k;\varepsilon_k)
+
R^2\varepsilon_{k+1} \\
&\leq
\frac12\varepsilon_k
+
R^2\varepsilon_{k+1}
=
(R^2+1)\varepsilon_{k+1}.
\end{aligned}
$$
This proves \eqref{eq:Ek-ek}.

It remains to relate $\varepsilon_k$ to $N_k$. Algorithm~\ref{alg: restart perturbed ATPD} gives
$$
m_i=C_0\varepsilon_i^{-1/2},
\qquad
C_0=
\left(\frac{4}{\theta}\sqrt{\frac{2L_f}{1-\delta}}+\sqrt{\varepsilon_0}\right)
\ln(4(1-\theta)^{-1}(R^2+1)).
$$
Thus
$$
N_k
=
C_0\sum_{i=0}^{k}\varepsilon_i^{-1/2}
=
\frac{\sqrt{2}C_0}{\sqrt{2}-1}
\left(\varepsilon_k^{-1/2}-\frac{\varepsilon_0^{-1/2}}{\sqrt{2}}\right),
$$
which gives \eqref{eq:rate}. Rounding adds at most one iteration per stage, so the complexity estimate follows. \qed
\end{proof}

\begin{remark}[\revision{Practical choice of homotopy parameters}]
In Algorithm~\ref{alg: restart perturbed ATPD}, the parameters $\varepsilon_0$ and $m_0$ are prescribed by the theory. In practice, we use the simpler choices
$$
\theta=1,\qquad\varepsilon_0=L_f,\qquad m_0\in[10,100].
$$
These choices define a practical variant with $m_0$ specified independently of the theoretical formula; they are outside the assumptions of Theorem~\ref{thm: global convergence rate}. This variant is robust in our experiments.
%The number $m_k$ is the maximal number of inner ATPD iterations at stage $k$. Since Algorithm~\ref{alg:ATPD} includes a stopping criterion for the inner solve, the inner loop may stop before completing all $m_k$ iterations. Thus $m_k$ should be interpreted as an upper bound rather than as a fixed iteration count.
\end{remark}

\subsection{An Accelerated Augmented Lagrangian Method}\label{sec: ALM}

We relate ATPD to a linearized augmented Lagrangian method. For a convex and $L_f$-smooth function $f$, 
% We also discuss the case where $f$ is convex on $\mathcal V$ but strongly convex only on the feasible set
% $$
% \mathcal K:=\{x\in\mathcal V:Bx=b\}.
% $$
consider the augmented Lagrangian problem
\begin{equation}\label{eq: augmented Lagrangian}
\min_{x\in\mathcal V}\max_{\lambda\in\mathcal Q}
\mathcal L_\beta(x,\lambda)
:=
f(x)
+\frac{\beta}{2}\|Bx-b\|^2
+\inprd{\lambda,Bx-b},
\end{equation}
where $\beta>0$. Define
$$
f_\beta(x):=f(x)+\frac{\beta}{2}\|Bx-b\|^2.
$$
Suppose $f_\beta$ is $\mu_{f_\beta}$-strongly convex. Its strong convexity and smoothness constants satisfy
% If $f$ is $\mu_{\mathcal K}$-strongly convex on $\mathcal K$, then
% $$
% 0<
% \min\{\mu_{\mathcal K},\,\beta\lambda_{\min}(BB^{\intercal})\}
% \leq
% \mu_{f_\beta}
% \leq
% L_{f_\beta}
% \leq
% L_f+\beta\lambda_{\max}(BB^{\intercal}).
% $$
$$
0<
\mu_{f_\beta}
\leq
L_{f_\beta}
\leq
L_f+\beta\lambda_{\max}(BB^{\intercal}).
$$
Applying ATPD~\eqref{eq:ATPD original} to \eqref{eq: augmented Lagrangian}, with $f_\beta$ in place of $f$, gives a linearized augmented Lagrangian method.

Theorem~\ref{thm: linear convergence of ATPD} gives the complexity estimate. With accelerated iterative solvers for the inner linear systems, the method needs
$$
\nabla f_\beta:\ 
\mathcal O\left(\sqrt{\kappa(f_\beta)}\log\frac{1}{\varepsilon}\right),
\qquad
B \text{ or } B^{\intercal}:\ 
\mathcal O\left(\sqrt{\kappa(f_\beta)\kappa(BB^{\intercal})}\log\frac{1}{\delta}\log\frac{1}{\varepsilon}\right)
$$
gradient evaluations and matrix--vector products, respectively.

Accelerated sublinear rates have been established for related ALM and ADMM methods~\cite{xu2017accelerated,ouyang2015accelerated,boct2023fast}. Here, ATPD gives an accelerated linear rate when $\mu_{f_{\beta}}>0$. When $\mu_{f_{\beta}}=0$, applying the homotopy ATPD method in Algorithm~\ref{alg: restart perturbed ATPD} to the augmented Lagrangian gives the accelerated $\mathcal O(1/N^2)$ decay of the modified Lyapunov energy after $N$ inner iterations under the uniform all-stage boundedness assumption of Theorem~\ref{thm: global convergence rate}.

%%%%%%%%%%%%%%%%%%%%%%%%%%%%
\section{Numerical Results}\label{sec: Numerical Results}
%%%%%%%%%%%%%%%%%%%%%%%%%%%%
In this section, we present numerical results for the proposed ATPD algorithms applied to partial differential equations, policy evaluation, and compressed sensing. 

\subsection{Stokes Equations}

We use the ill-conditioned Stokes equations as a benchmark to demonstrate the acceleration of the ATPD method. While multigrid methods for saddle-point systems can yield optimal solvers \cite{WangChen2013MGStokesLSC}, we focus here on the performance of our first-order methods and show the acceleration.

Consider the steady Stokes equations
\begin{equation}\label{eq: stokes}
\begin{aligned}
-\Delta \mathbf{u} + \nabla p &= \mathbf{f}, \\
\nabla \cdot \mathbf{u} &= 0,
\end{aligned}
\end{equation}
on the unit square $\Omega=(0,1)^2$ with homogeneous Dirichlet boundary conditions $\mathbf{u}|_{\partial \Omega} = 0$. Here $\mathbf{u}$ is the velocity and $p$ is the pressure. Let $H_0^1(\Omega)$ be the standard Sobolev space and $L_0^2(\Omega):=\{q\in L^2(\Omega): \int_{\Omega} q = 0\}$. We use a uniform triangulation $\mathcal T_h$ of $\Omega$ with mesh size $h$, and Taylor--Hood ($P_2$--$P_1$) finite element spaces $\mathcal V_h\subset [H_0^1(\Omega)]^2$ and $\mathcal Q_h\subset L_0^2(\Omega)$.

The discrete Stokes system is the KKT condition for the following constrained minimization problem:
\begin{equation}\label{eq:Stokes opt}
\begin{aligned}
\min_{\mathbf{u}_h \in \mathcal V_h}\quad & f_h(\mathbf{u}_h)
:= \frac{1}{2}\int_{\Omega} |\nabla \mathbf{u}_h|^2 \dd x - \inprd{\mathbf f, \mathbf u_h} \\
\text{subject to}\quad & (\nabla\cdot \mathbf u_h, q_h) = 0, \quad \forall~q_h \in \mathcal Q_h.
\end{aligned}
\end{equation}
The pressure $p$ in \eqref{eq: stokes} acts as the Lagrange multiplier for the divergence-free constraint. In the Euclidean norm of the coefficient vector in the standard finite element basis, it is well known that $L_{f_h}=\mathcal O(1)$ is independent of $h$, while $\mu_{f_h}=\mathcal O(h^2)$, giving a condition number $\kappa(f_h) = C h^{-2}$. Let $B_h$ be the matrix representation of the discrete divergence operator. Geometric multigrid provides an efficient preconditioner $\IQ$ for solving $B_h B_h^{\intercal}$.

We implement the problem using the software package {\it i}FEM~\cite{chen2008ifem} and terminate the algorithms when the relative residual falls below $10^{-8}$. We estimate $\mu_{f_h}$ and $L_{f_h}$ from a coarse mesh $H=1/4$ using
$$ \mu_{f_h} = \mu_{f_H}(h/H)^2, \quad L_{f_h} = L_{f_H}, $$
where $\mu_{f_H}$ and $L_{f_H}$ are found by direct factorization.

We compare Algorithm~\ref{alg: TPD} (TPD) and Algorithm~\ref{alg:ATPD} (ATPD). In each iteration, we apply one V-cycle multigrid as $\IQ^{-1}$, which costs $\mathcal O(N)$ flops for $N$ degrees of freedom. We set the step size $\alpha=1.5$ for TPD and $\alpha=1/\sqrt{\kappa(f_h)}$ for ATPD. These choices are larger than the conservative theoretical bounds, but both methods converge and the solutions reach the expected accuracy.

\begin{figure}[!htb]
\begin{minipage}{0.62\textwidth}
\centering
\renewcommand{\arraystretch}{1.25}
\resizebox{\textwidth}{!}{%
\begin{tabular}{@{} c c  c c  c c @{}}
\toprule
& & \multicolumn{2}{c}{TPD Algorithm~\ref{alg: TPD}}
& \multicolumn{2}{c}{ATPD Algorithm~\ref{alg:ATPD}} \\
$h$ & DoF ($N$) & Iterations & Time & Iterations & Time \\
\midrule
$1/16$  & 2,467   & 2,065   & 0.42s   & 1,112   & 0.22s \\
$1/32$  & 9,539   & 8,262   & 3.87s   & 1,745   & 0.83s \\
$1/64$  & 37,507  & 33,045  & 58.81s  & 3,647   & 6.17s \\
$1/128$ & 148,739 & 132,159 & 1,176s  & 7,079   & 54.76s \\
\bottomrule
\end{tabular}
}
\captionof{table}{Performance comparison of the TPD and ATPD algorithms for solving the Stokes equations.} 
%The iteration count of the TPD algorithm scales as $\mathcal O(\kappa(f_h))=\mathcal O(N)$, whereas the ATPD algorithm exhibits accelerated convergence with iteration count scaling as $\mathcal O(\sqrt{\kappa(f_h)})$.}
\label{table:stokes-results}
\end{minipage}
\hfill
\begin{minipage}{0.35\textwidth}
\centering
 \includegraphics[width = \linewidth]{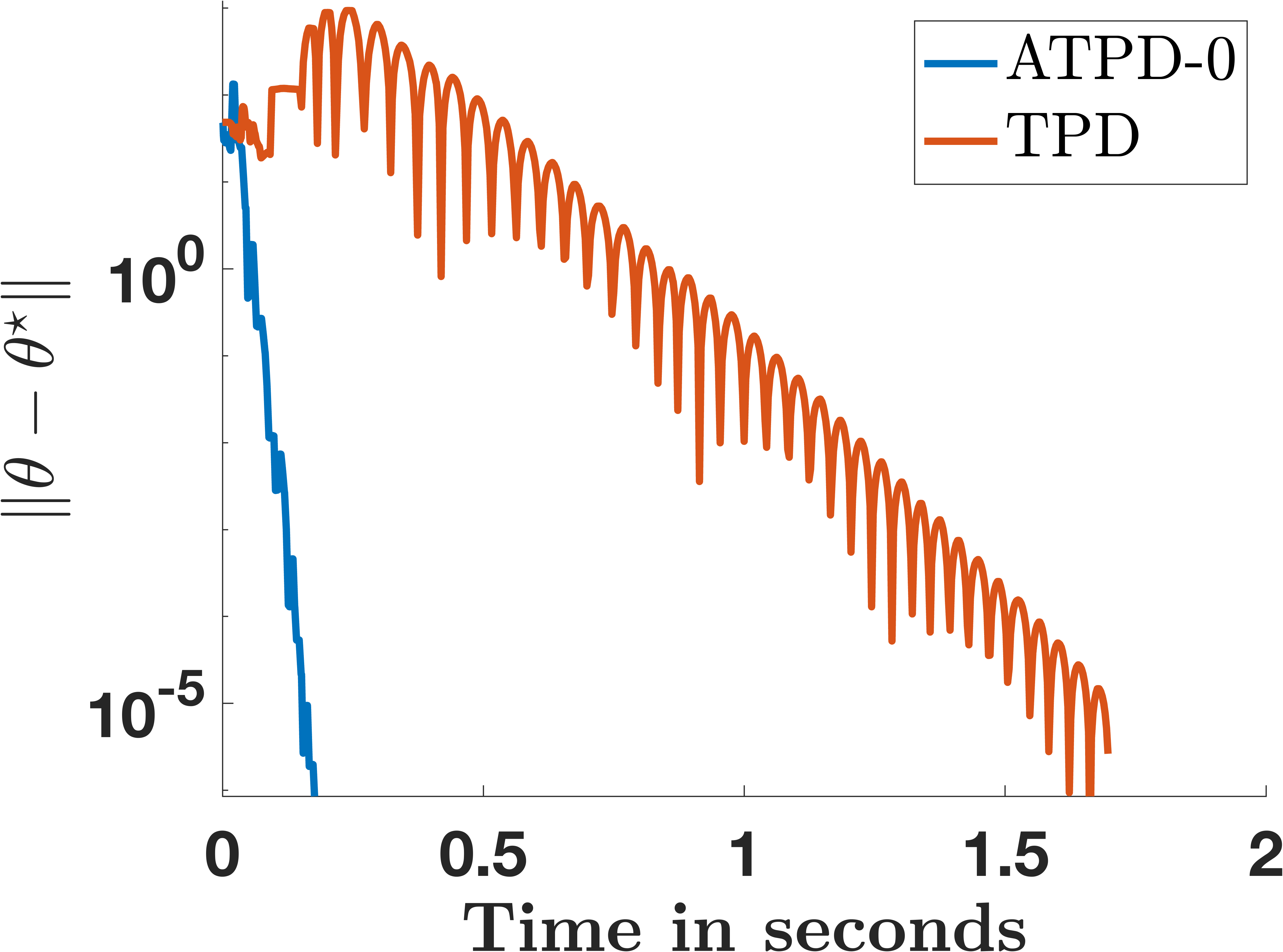}
 \caption{TPD and ATPD algorithms for solving the policy evaluation problem.}
  \label{fig:ATPD-policy evaluation time}
\end{minipage}
\end{figure}

Table~\ref{table:stokes-results} lists iteration counts and CPU times for various $h$ and degrees of freedom $N:=\operatorname{dim}(\mV_h) + \operatorname{dim}(\mQ_h)$. The TPD iteration count grows as $\mathcal O(\kappa(f_h)) = \mathcal O(h^{-2}) = \mathcal O(N)$. By contrast, ATPD shows accelerated growth of $\mathcal O(\sqrt{\kappa(f_h)}) = \mathcal O(h^{-1}) = \mathcal O(\sqrt{N})$. With an inner cost of $\mathcal O(N)$ per outer iteration, the total complexity is $\mathcal O(N^2)$ for TPD and $\mathcal O(N^{3/2})$ for ATPD.

\vspace{-0.5cm}

\subsection{Policy Evaluation}
We consider the policy evaluation problem in reinforcement learning~\cite{du2017stochastic}. For a Markov decision process (MDP) under policy $\pi$ with feature map $\phi: \mathcal S \to \mathbb{R}^d$, minimizing the mean squared projected Bellman error (MSPBE) with a linear value function approximation is a weighted least-squares problem~\cite{dann2014policy}:
\begin{equation}\label{eq:MSPBE}
    \min_{\theta} \text{MSPBE}(\theta) = \frac{1}{2}\|A\theta - b\|^2_{C^{-1}},
\end{equation}
where
$$
A = \mathbb{E}[\phi(s_t)(\phi(s_t) - \gamma \phi(s_{t+1}))^{\intercal}], \quad b = \mathbb{E}[\phi(s_t)r_t], \quad C = \mathbb{E}[\phi(s_t)\phi(s_t)^{\intercal}].
$$
Here $\gamma \in [0, 1)$ is the discount factor and expectations are taken with respect to the stationary distribution. We assume $C\succ0$ and $A$ is nonsingular. By introducing an auxiliary variable $z$ such that $Cz = b-A\theta$, the first-order optimality condition for the MSPBE problem \eqref{eq:MSPBE} is equivalent to solving the linear system:
\begin{equation}\label{eq:MSPBElinear}
    \begin{pmatrix}
C & A \\
A^{\intercal} & 0
\end{pmatrix}
\begin{pmatrix}
z \\
\theta
\end{pmatrix} = 
\begin{pmatrix}
b \\
0
\end{pmatrix}.
\end{equation}
\revision{In particular, \eqref{eq:MSPBElinear} is the KKT system of the problem in the form of \eqref{eq: affine equality constrained optimization system}:
\begin{equation}\label{eq:MSPBEconstrained}
    \min f(z):= \frac{1}{2}\|z\|^2_C - b^{\intercal }z, \quad \text{ s.t. } A^{\intercal }z = 0.
\end{equation}
}

We test an MDP problem with $400$ states, $10$ actions, and a $201$-dimensional feature vector. The transition probabilities, rewards, and features are sampled from uniform distributions, and the discount factor is $\gamma=0.95$.

Since $\mu_f$ is hard to estimate, we use the homotopy ATPD method, Algorithm~\ref{alg: restart perturbed ATPD}, denoted by ATPD-0, with algorithm parameters $\theta = 1$, $\varepsilon_0=L_f$ and $m_0=80$. We estimate $L_f$ by applying the power method to $C$. Here $B=A^{\intercal}$, so the inner system has matrix $BB^{\intercal}=A^{\intercal}A$. The inner solver $\IQ^{-1}$ is PCG with a diagonal preconditioner. We stop the inner iteration when the relative residual reaches $\delta=10^{-1}$.

Figure~\ref{fig:ATPD-policy evaluation time} shows error convergence versus CPU time. ATPD-$0$ reaches small errors substantially faster than TPD and reduces error oscillations.

\subsection{Compressed Sensing}

We consider a compressed sensing experiment with regularized objective~\cite{salim2022optimal}:
\begin{equation}\label{eq: smooth compressed sensing}
    \min_{Ax = b} f(x) := \sum_{i=1}^n g\bigl(x_i\bigr), 
    \quad 
    g(x) = \sqrt{x^2+c^2}+\frac{c}{2}x^2,
\end{equation}
where $A \in \mathbb{R}^{m \times n}$, $b \in \mathbb{R}^m$, $m<n$, and $f$ is a $\mu_f$-strongly convex and $L_f$-smooth approximation of the $\ell_1$ norm. The goal is to recover a sparse vector $x^{\sharp}\in\mathbb{R}^n$ with $50$ random nonzero entries, all equal to $1$, from measurements $b=Ax^\sharp$.

Direct calculation gives
$$
\mu_f = c, 
\quad 
L_f = \frac{1}{c}+c,
\quad 
\kappa(f)=\frac{L_f}{\mu_f}=1+\frac{1}{c^2}.
$$
The parameter $c$ controls the smoothing. Smaller $c$ gives a closer approximation to the $\ell_1$ norm, but also gives a larger condition number.

Let $x^\star$ be the solution of~\eqref{eq: smooth compressed sensing}. If $f$ were exactly the $\ell_1$ norm, then $x^\star=x^\sharp$ under suitable sparse-recovery conditions on $A$. In our tests, we choose parameters for which $x^\star\neq x^\sharp$ to examine performance in an ill-conditioned smooth problem. We keep the two vectors separate.

\subsubsection{Inexact Inner Solve}

We first test synthetic data with $m=250$ and $n=1000$. We draw a matrix with independent Gaussian entries and then obtain $A$ by rescaling its nonzero singular values to lie in $[1/\sqrt{\chi},1]$ for a prescribed $\chi$.

We examine the effect of the inner-solver tolerance $\delta$. We set $\chi=10^4$ and $\kappa(f)=10^3$, and use the step size $\alpha=1/\sqrt{2\kappa(f)}$. The algorithm stops when the relative residual of~\eqref{eq: smooth compressed sensing} is below $10^{-8}$. The operator $\IQ^{-1}$ is applied by PCG with a diagonal preconditioner for $AA^{\intercal}$. The PCG iteration stops when its relative residual reaches $\delta\in\{0.5,10^{-1},10^{-2},10^{-4}\}$.

\begin{figure}
\begin{minipage}{0.32\textwidth}
\centering
 \includegraphics[width = 0.9\linewidth]{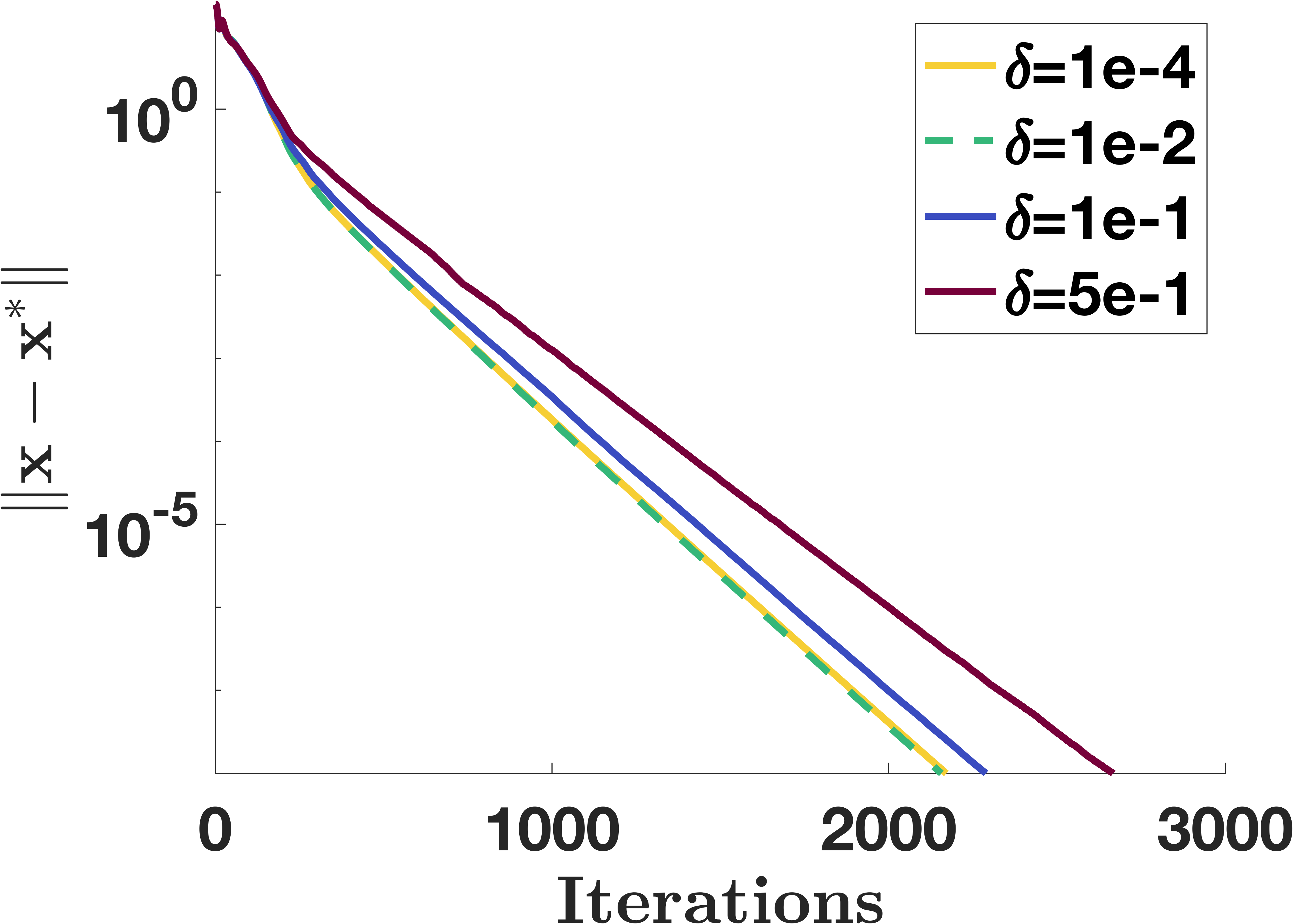}
\end{minipage}
\hfill
\begin{minipage}{0.32\textwidth}
\centering
 \includegraphics[width = 0.9\linewidth]{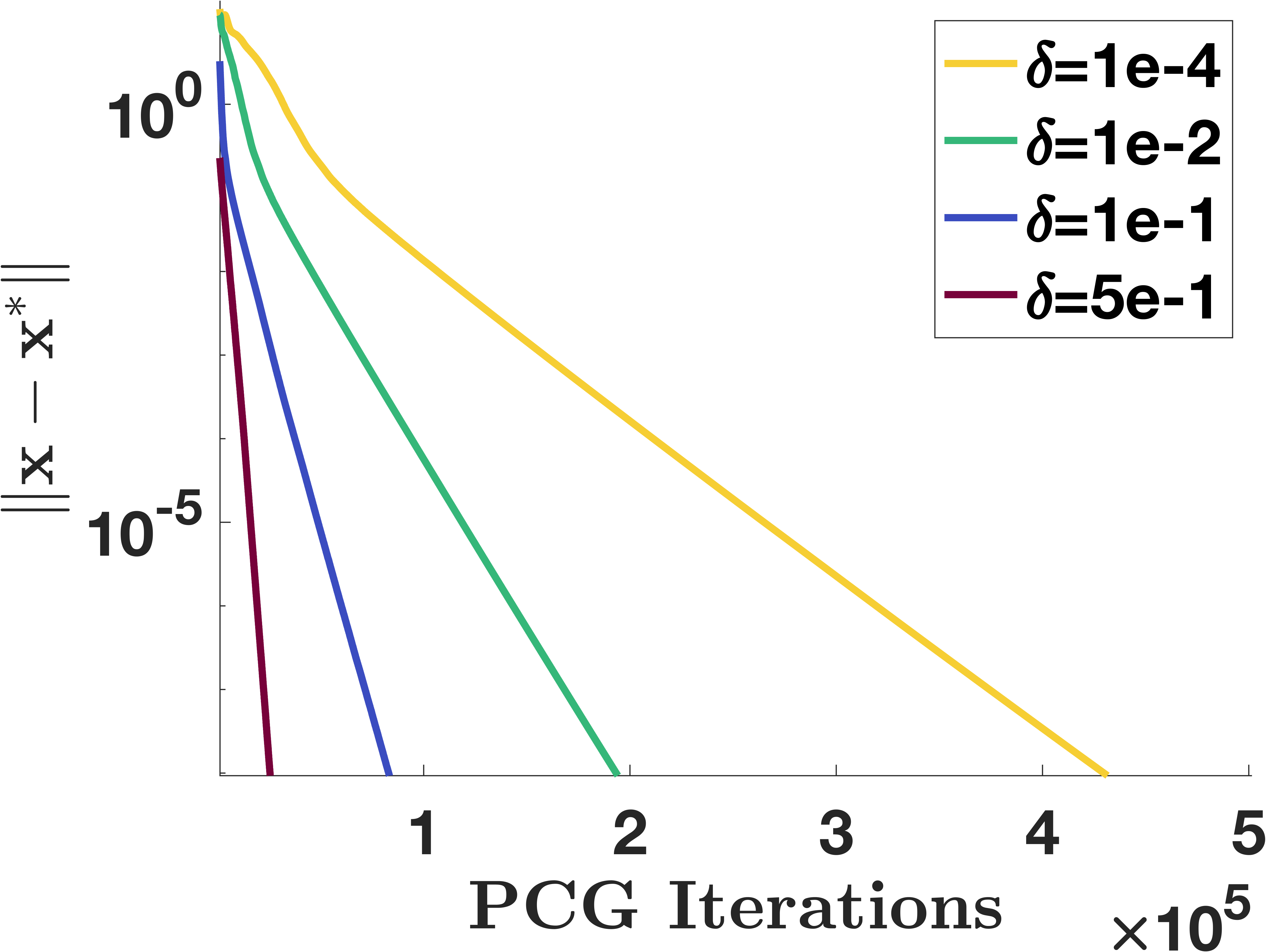}
\end{minipage}
\hfill
\begin{minipage}{0.32\textwidth}
\centering
 \includegraphics[width = 0.9\linewidth]{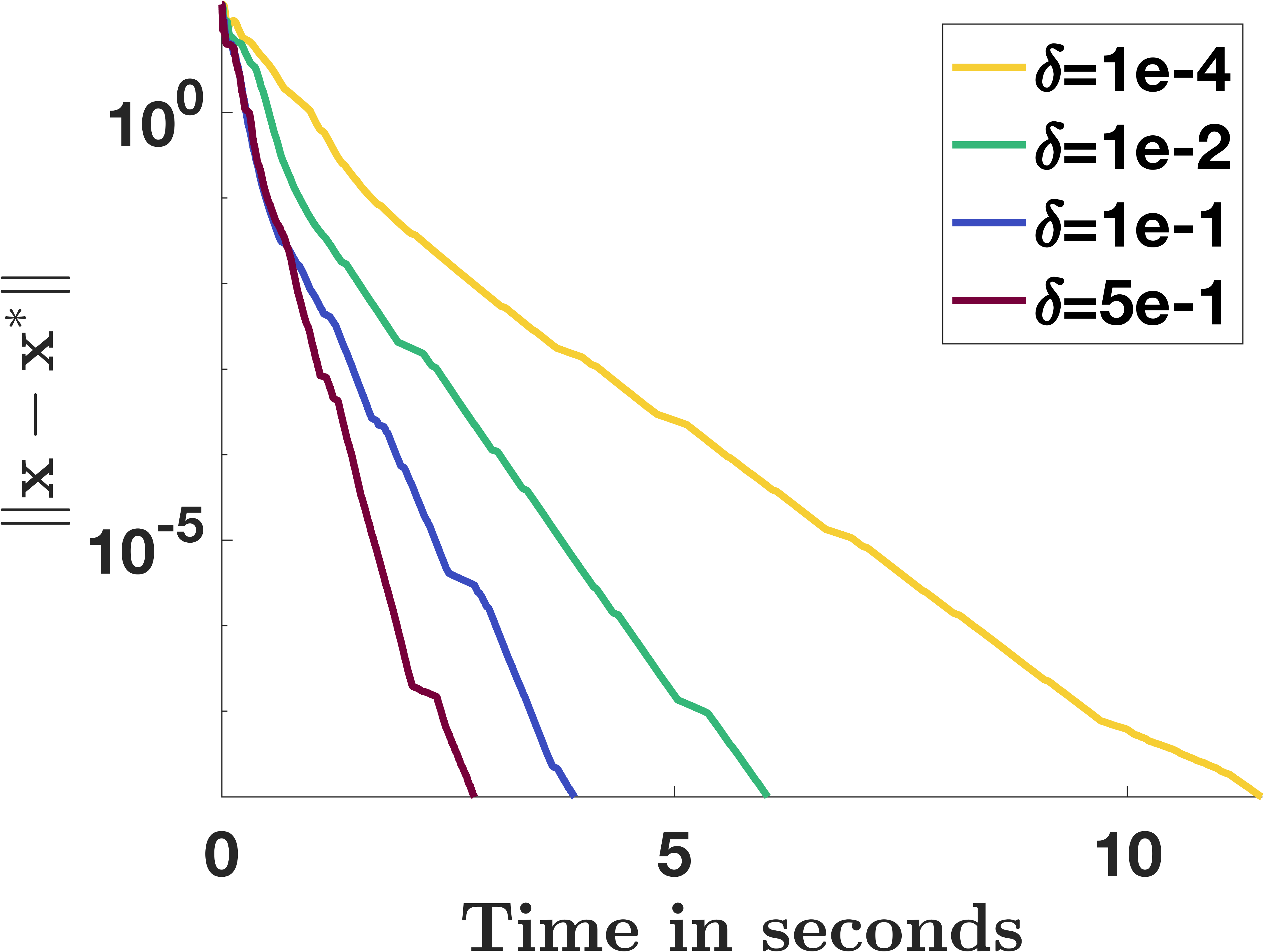}
\end{minipage}
\caption{Numerical results for the ATPD method solving the compressed sensing problem with various $\delta$ for the inner solve.}
 \label{fig:ATPD-compress delta comparison}
\end{figure}

Figure~\ref{fig:ATPD-compress delta comparison} shows the error decay, total PCG iterations, and CPU time. Smaller $\delta$ gives more accurate inner solves and thus fewer outer iterations. The observed outer iteration count changes only mildly with $\delta$. In contrast, the total number of PCG iterations and the runtime decrease substantially as $\delta$ increases. Looser inner tolerances substantially reduce the total PCG iterations in these experiments.

Since the inner solver is the main cost, a coarse tolerance such as $\delta=0.5$ is effective. For very small $\delta$, the method approaches accelerated projected gradient behavior. ATPD remains effective with coarse inner solves, which is its main advantage over projection-based methods.

\subsubsection{Comparison of ATPD Algorithms}

Next, we compare several variants of ATPD. We observe that using a value of $\mu_f$ larger than the theoretical lower bound can further accelerate convergence. We therefore tune this parameter and denote the method by ATPD-$\mu^\ast$. Since $\mu_f$ is usually unknown in practice, we also test the homotopy strategy in Algorithm~\ref{alg: restart perturbed ATPD}, denoted by ATPD-$0$. For reference, we include \revision{the algorithm (OptAlg) proposed in~\cite{borodich2025near}, the accelerated linearized augmented Lagrangian
method (ALALM) in \cite{he2026accelerated}} and the accelerated PAPC (APAPC) algorithm from~\cite{salim2022optimal}. OptAlg is near-optimal for a broader class of optimization problems, while APAPC is the first optimal algorithm for strongly convex minimization under affine constraints.

\begin{figure}
 \begin{minipage}{0.32\textwidth}
\centering
 \includegraphics[width = 0.9\linewidth]{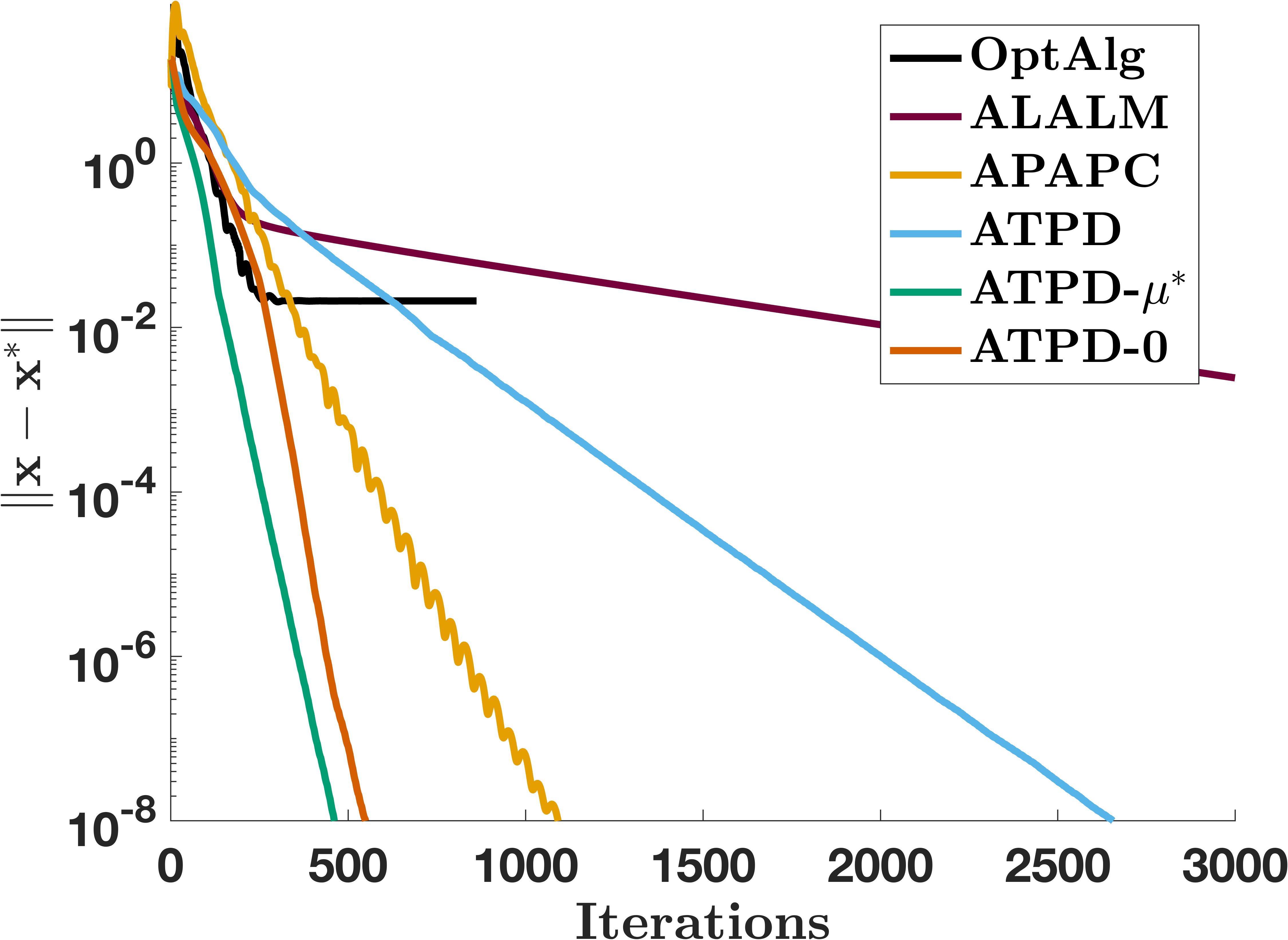}
\end{minipage}
\hfill
\begin{minipage}{0.32\textwidth}
\centering
 \includegraphics[width = 0.9\linewidth]{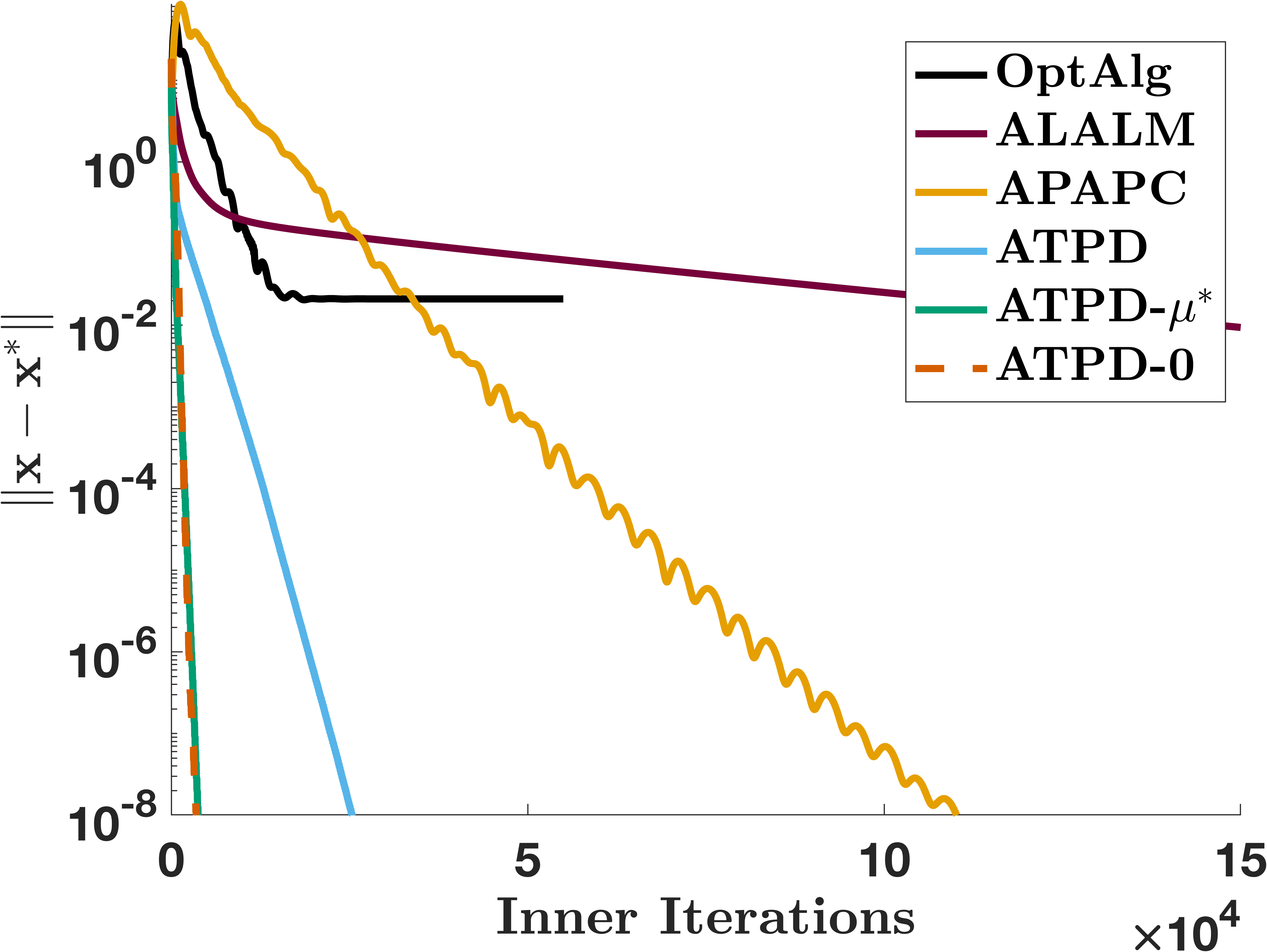}
\end{minipage}
\hfill
\begin{minipage}{0.32\textwidth}
\centering
 \includegraphics[width = 0.9\linewidth]{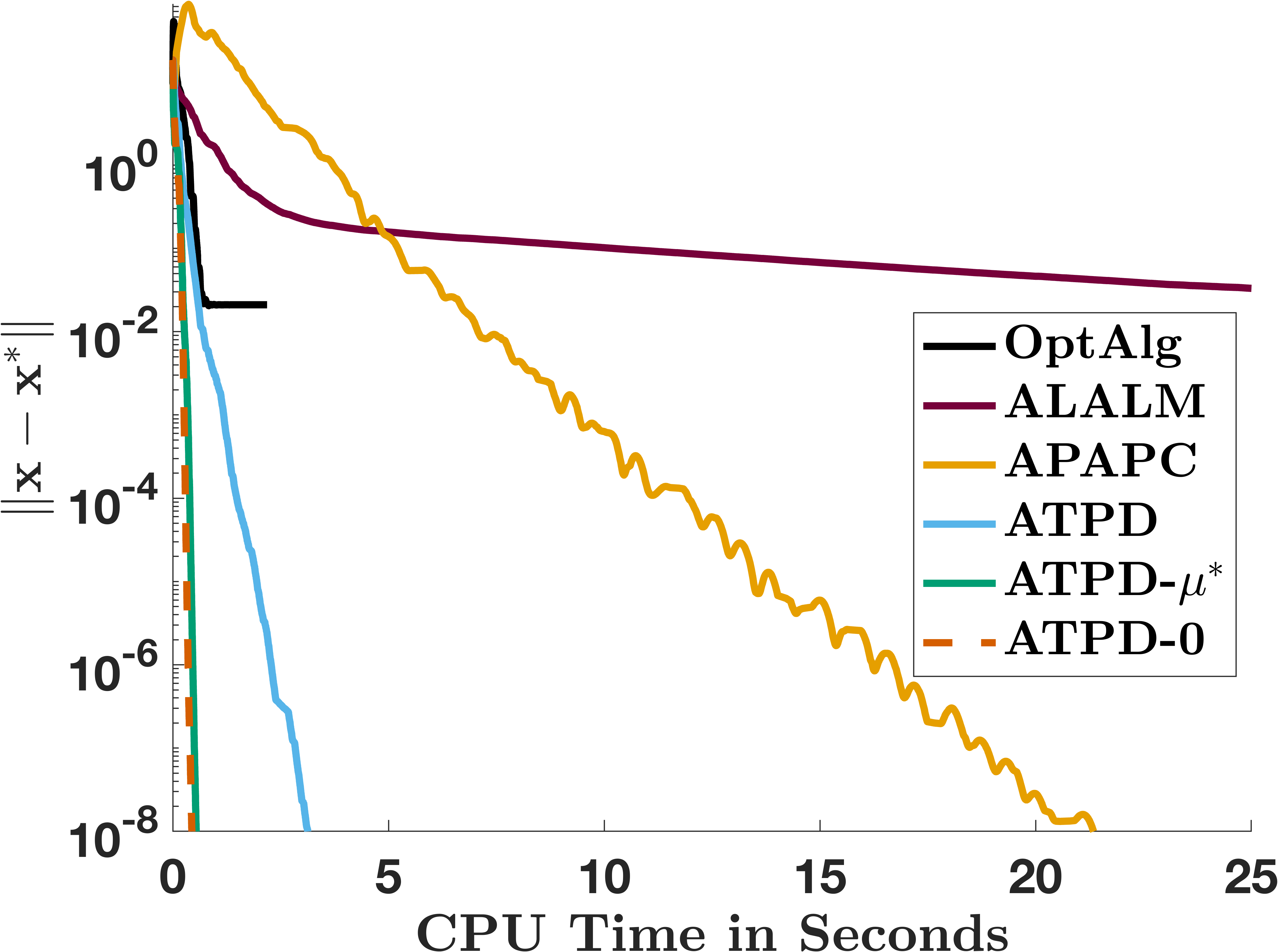}
\end{minipage}
\caption{Numerical results for variants of the ATPD method to solve the compressed sensing problem. OptAlg: the accelerated algorithm proposed in \cite{borodich2025near}; {\color{algdarkred} ALALM:} the accelerated linearized augmented Lagrangian
method in \cite{he2026accelerated}; {\color{algorange}  APAPC:} the accelerated PAPC algorithm in \cite{salim2022optimal}; {\color{algskyblue} ATPD:} Algorithm \ref{alg:ATPD}; {\color{alggreen} ATPD-${\mu}^*$:} Algorithm \ref{alg:ATPD} with tuned $\mu_f$; {\color{algvermillion} ATPD-$0$:} Algorithm \ref{alg: restart perturbed ATPD}.}
 \label{fig:ATPD-compress variants comparison}
\end{figure}

Figure~\ref{fig:ATPD-compress variants comparison} reports the outer iteration count, total inner iterations, and CPU time. The ATPD variants and APAPC exhibit approximately linear error decay over the plotted range; OptAlg reaches an error plateau, while ALALM decreases more slowly. \revision{OptAlg solves a perturbed regularized problem rather than \eqref{eq: smooth compressed sensing} directly, so its limit differs slightly from the true solution.  ALALM converges more slowly after hundreds of iterations as the penalty parameter becomes very large.} APAPC uses fewer outer iterations than ATPD, but its total runtime is much larger because it needs many inner-loop iterations to attain its optimal rate. In contrast, ATPD gives comparable or better performance with far fewer inner iterations.

Both ATPD-$\mu^\ast$ and ATPD-$0$ further improve the baseline ATPD method. For ATPD-$\mu^\ast$, setting $\mu^\ast = 50\,\mu_f$ improves convergence in this experiment, although this parameter choice is outside the stated convergence guarantee. For ATPD-$0$, we set $\theta =1$, $\varepsilon_0 = L_f/3$ and choose $m_0=80$. ATPD-$0$ avoids estimating $\mu_f$ by using the homotopy parameter $\varepsilon$ while maintaining robust performance.

\subsubsection{Benchmark Signal Reconstruction}
\begin{figure}[htbp]
\centering
    % --- (a) Sparco 9: Low resolution ---
    \begin{minipage}{0.32\textwidth}
        \centering
        \includegraphics[width=0.95\linewidth]{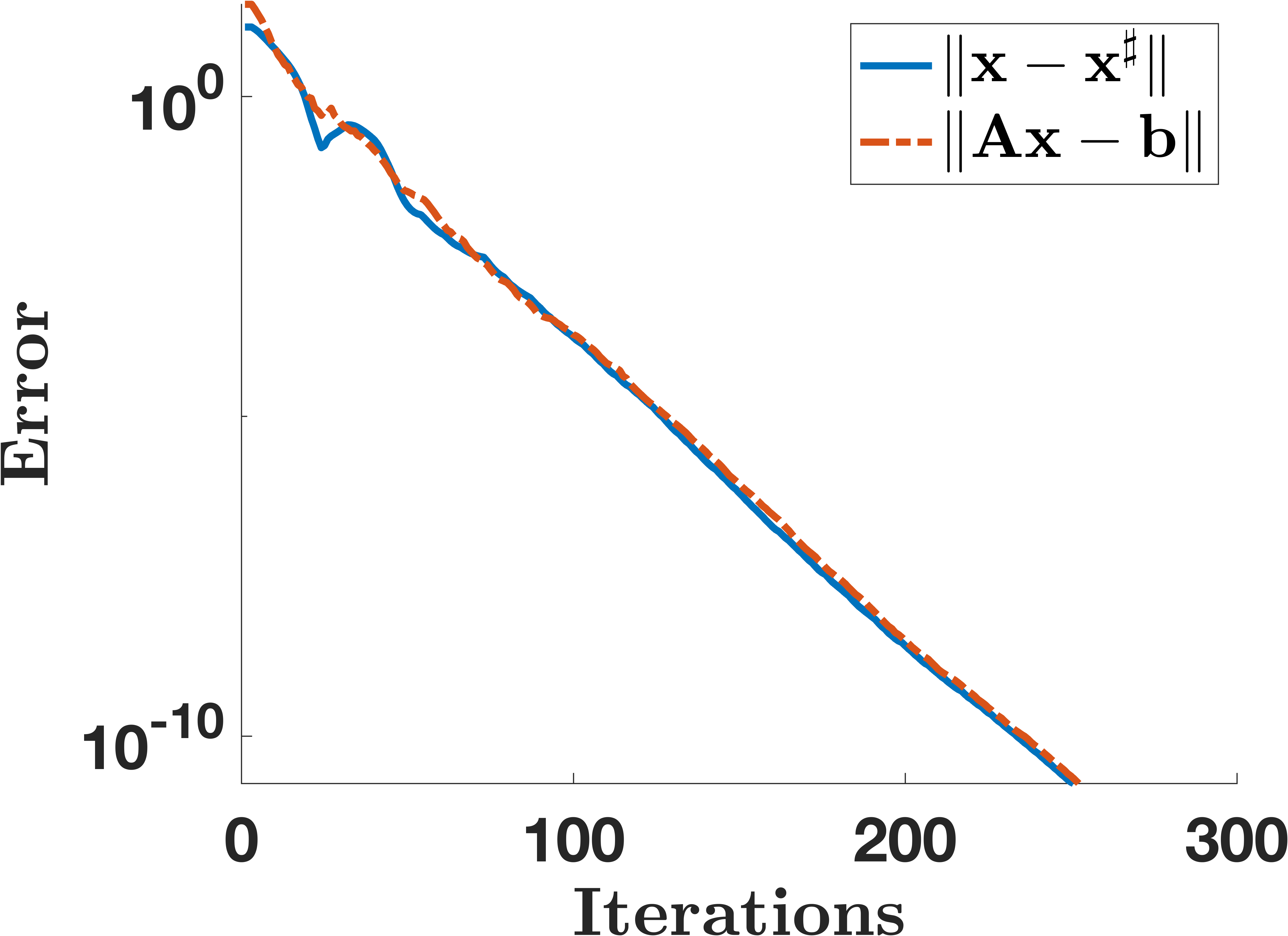}
    \end{minipage}
    \hfill
    \begin{minipage}{0.32\textwidth}
        \centering
        \includegraphics[width=0.85\linewidth]{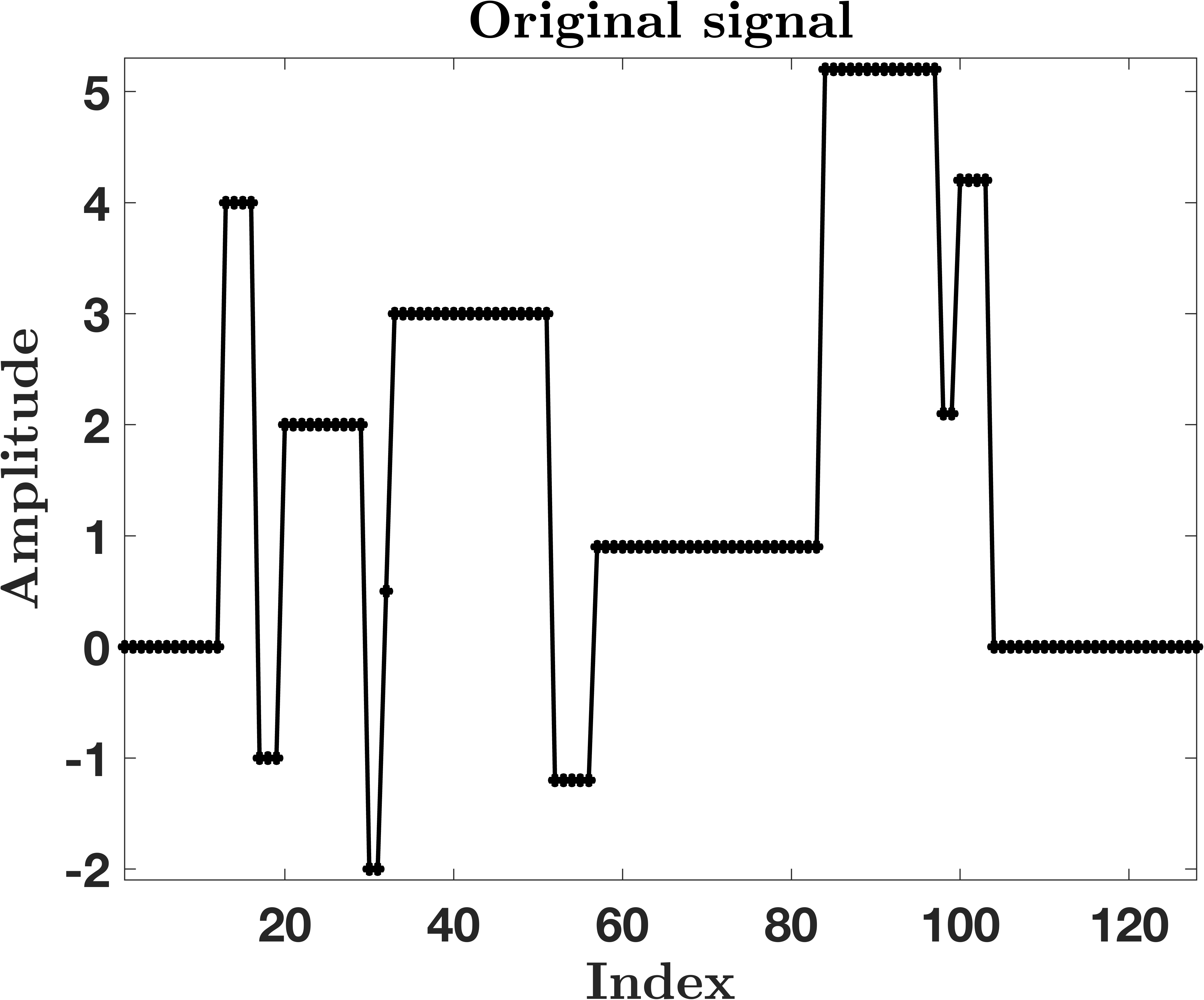}
    \end{minipage}
    \hfill
    \begin{minipage}{0.32\textwidth}
        \centering
        \includegraphics[width=0.85\linewidth]{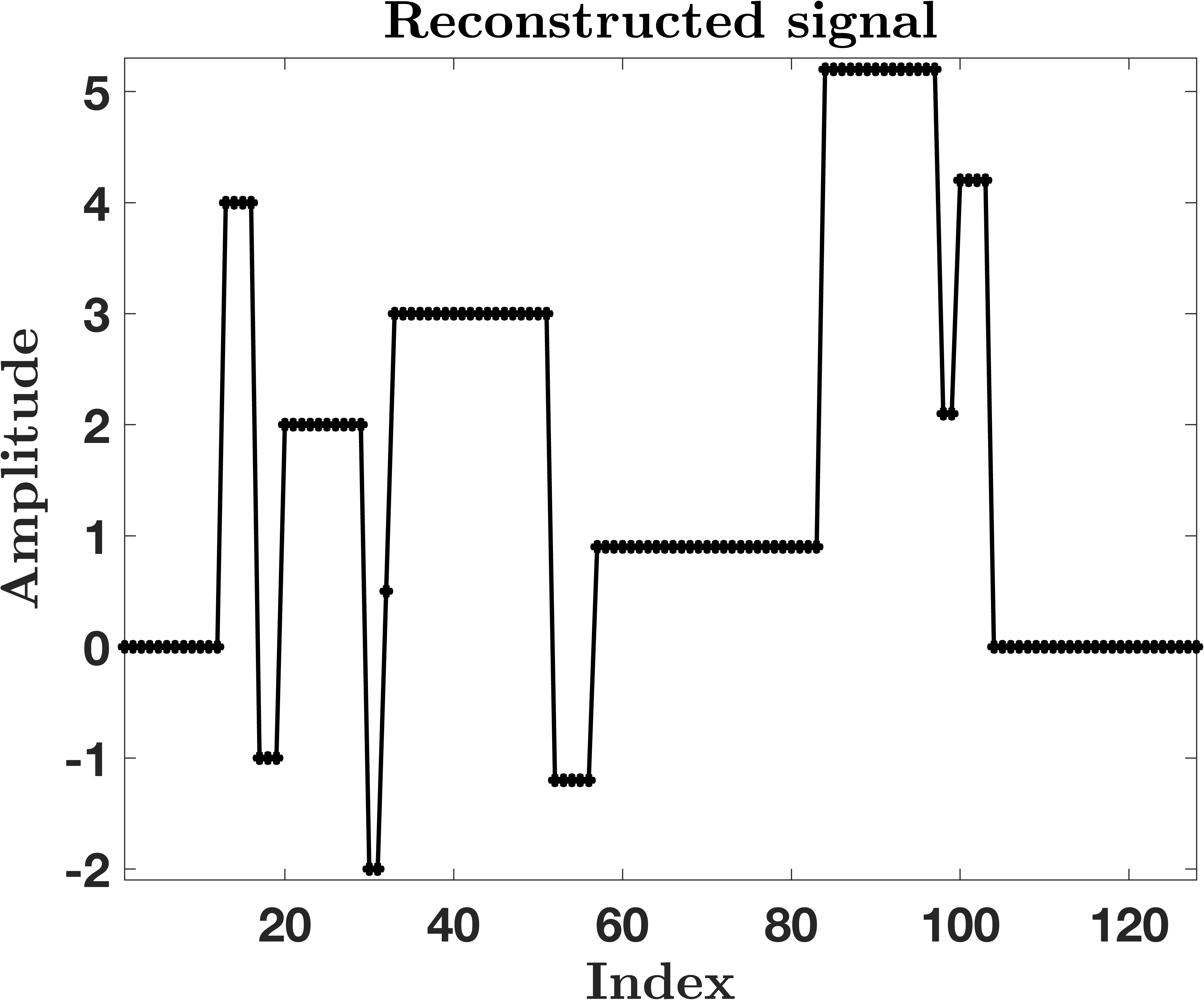}
    \end{minipage}
    
    \medskip % Vertical space between rows
    
    % --- (b) Sparco 10: High resolution ---
    \begin{minipage}{0.32\textwidth}
        \centering
        \includegraphics[width=0.95\linewidth]{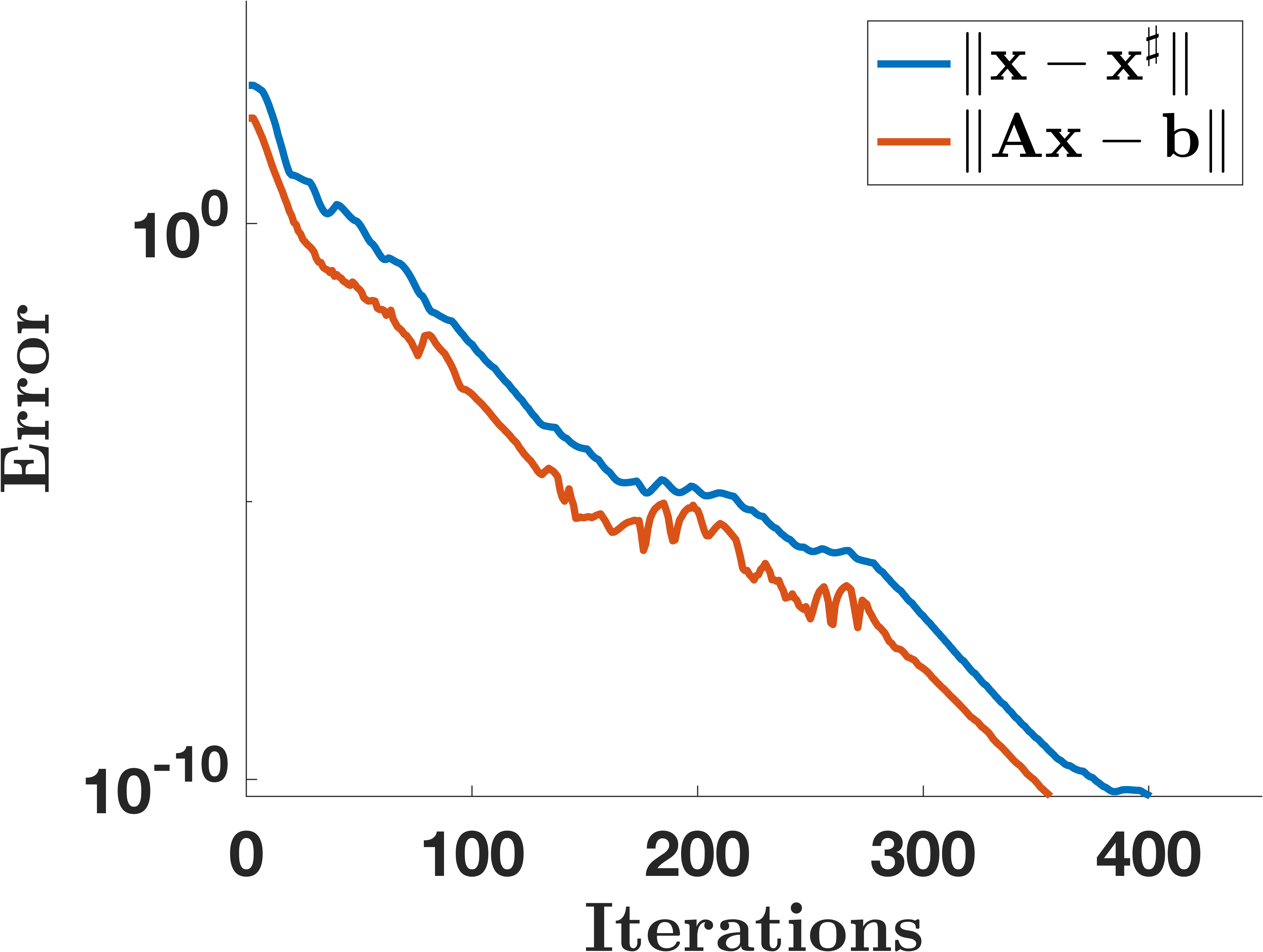}
    \end{minipage}
    \hfill
    \begin{minipage}{0.32\textwidth}
        \centering
        \includegraphics[width=0.85\linewidth]{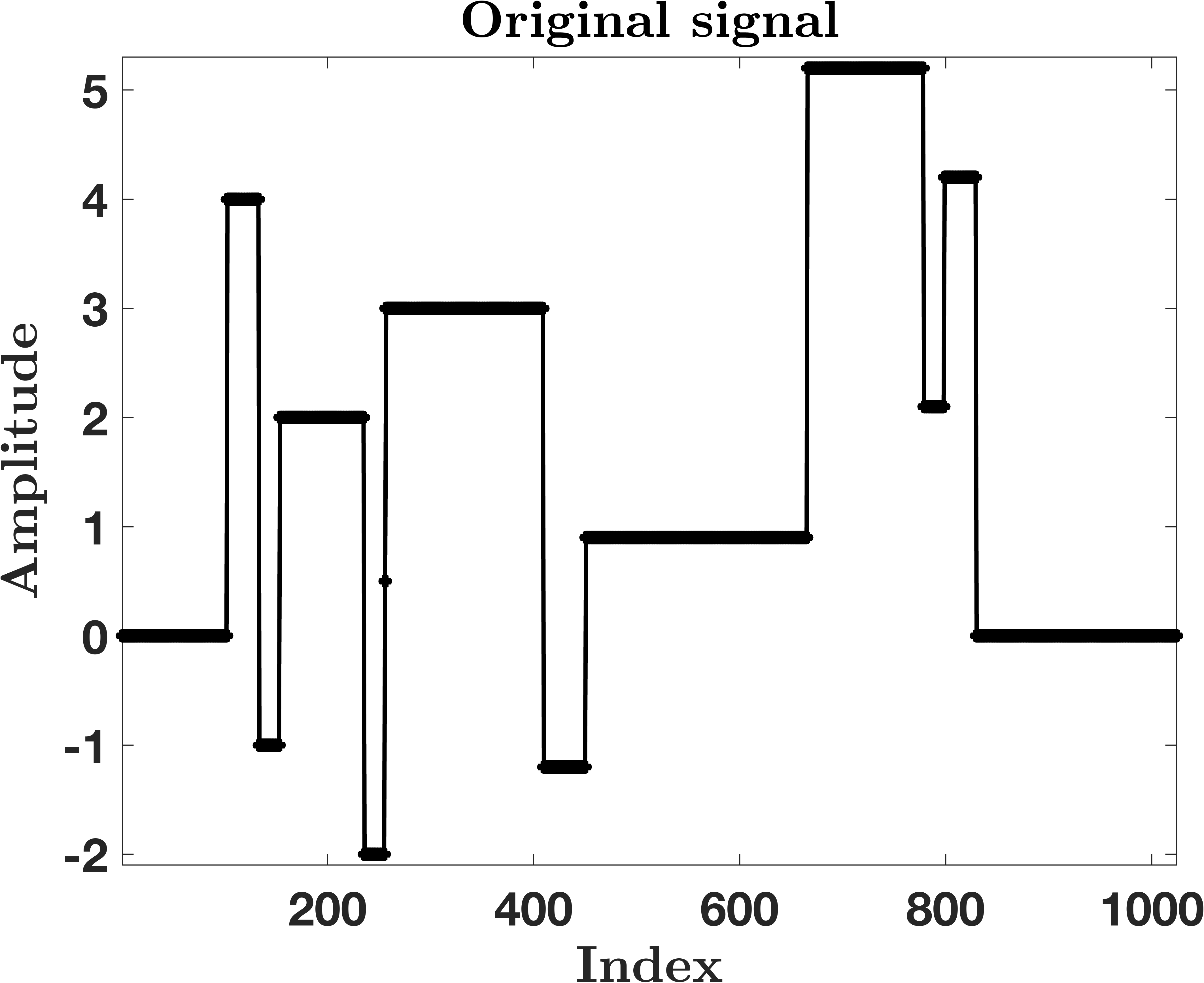}
    \end{minipage}
    \hfill
    \begin{minipage}{0.32\textwidth}
        \centering
        \includegraphics[width=0.85\linewidth]{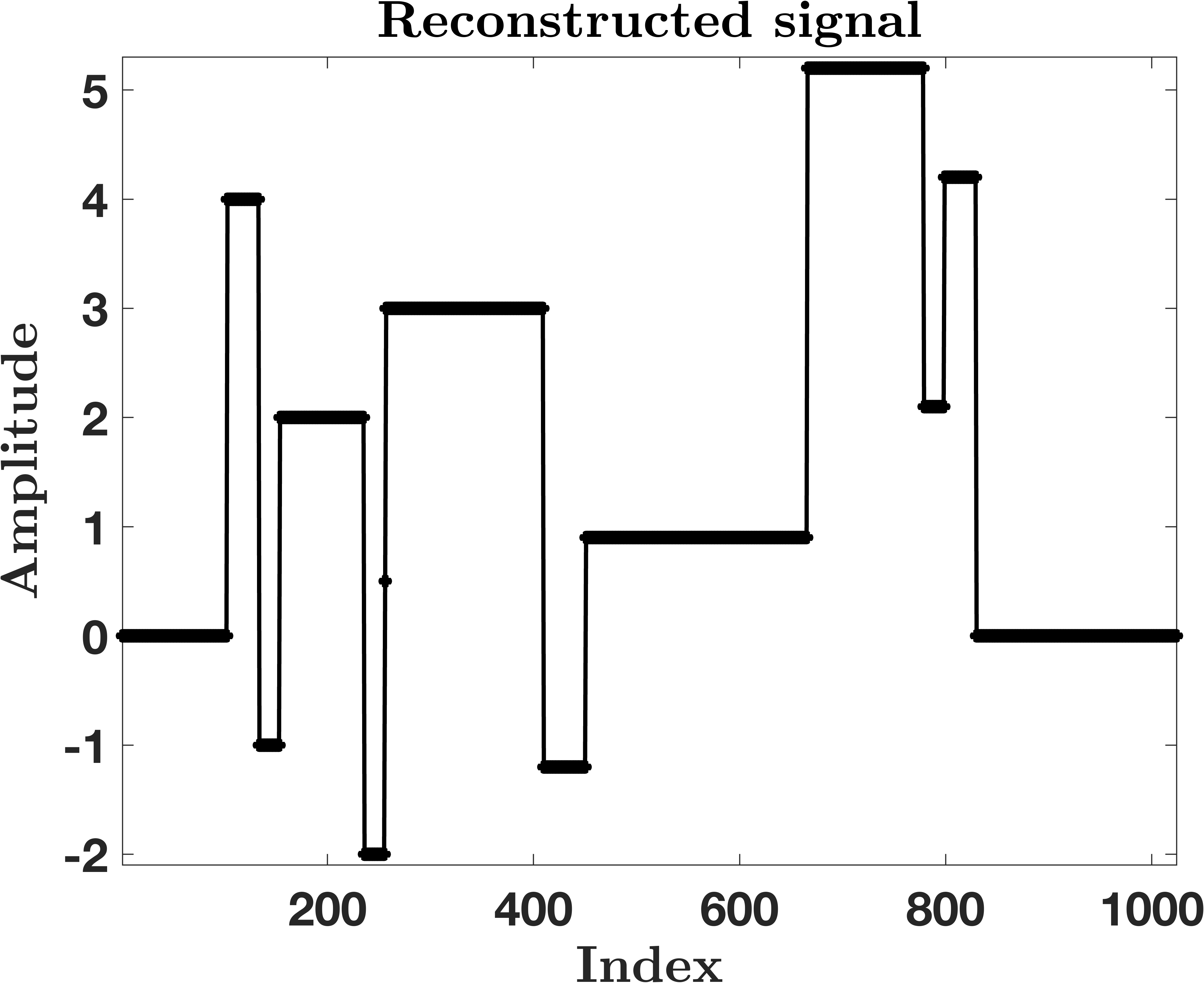}
    \end{minipage}

    \medskip % Vertical space between rows

    % --- (c) Sparco 11: High accuracy ---
 \begin{minipage}{0.32\textwidth}
\centering
 \includegraphics[width = 0.95\linewidth]{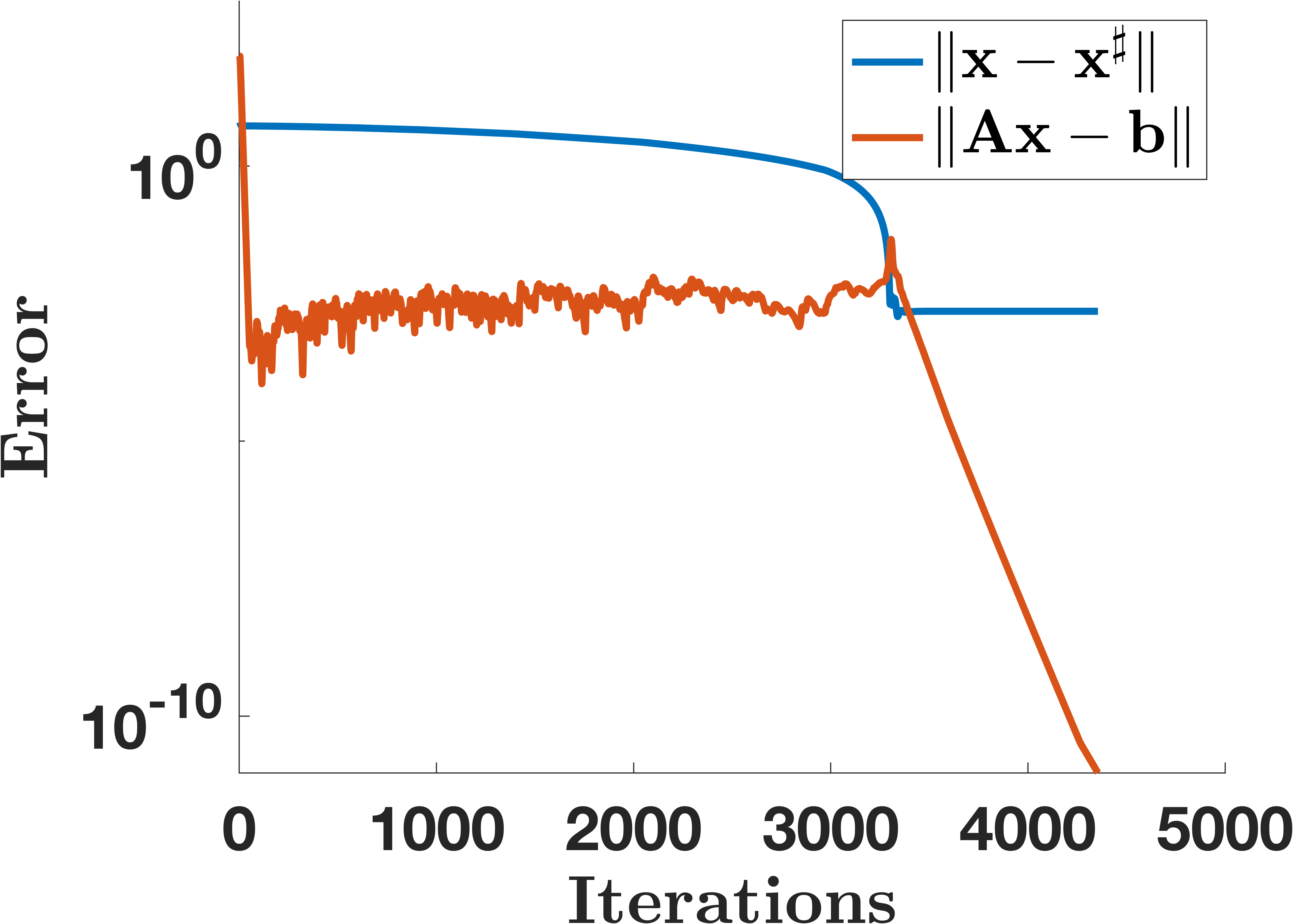}
    % \caption{Error v.s. Iteration steps. }
\end{minipage}
\hfill
\begin{minipage}{0.32\textwidth}
\centering
 \includegraphics[width = 0.85\linewidth]{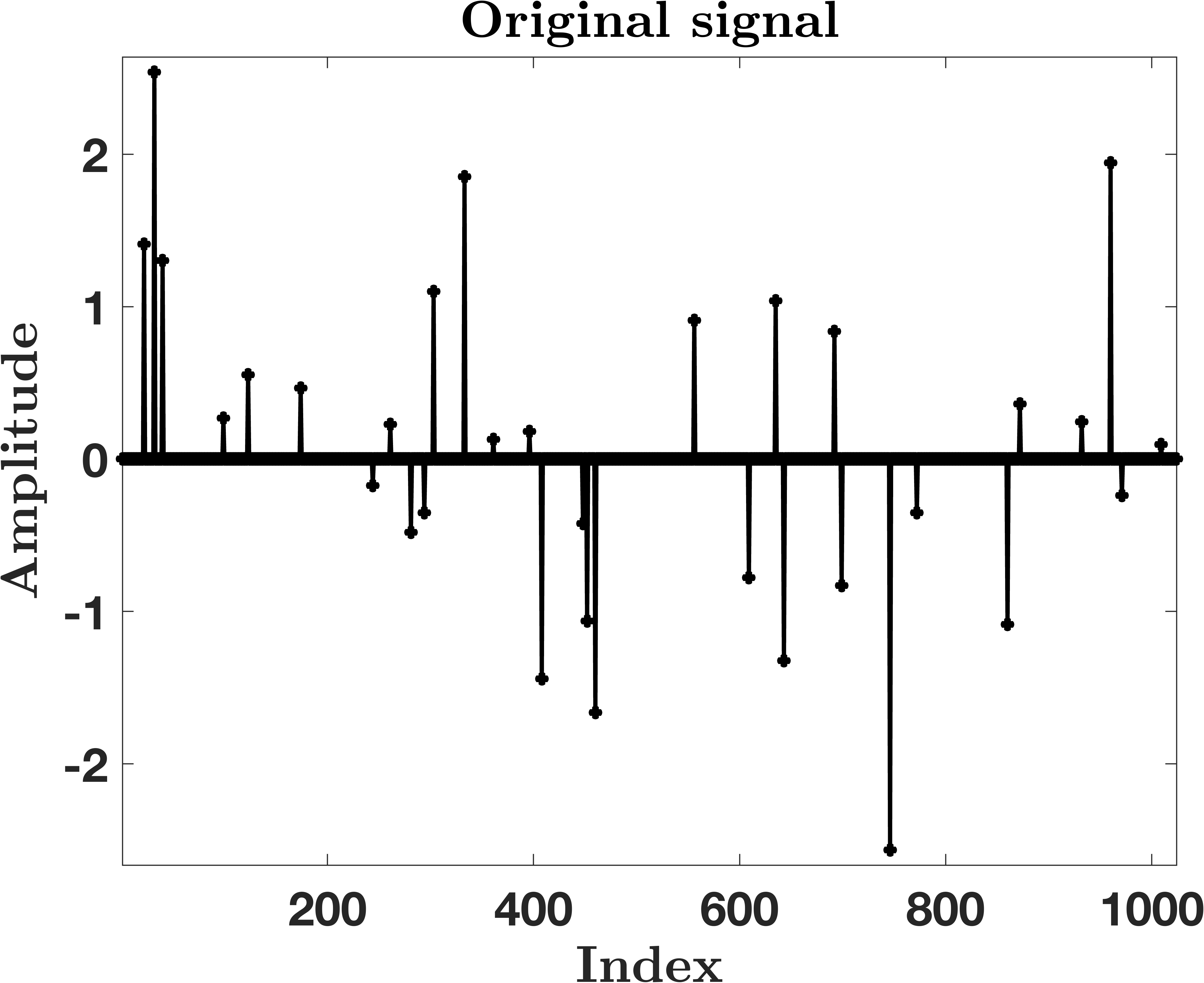}
    % \caption{Error v.s. Inner iterations. }
\end{minipage}
\hfill
\begin{minipage}{0.32\textwidth}
\centering
 \includegraphics[width = 0.85\linewidth]{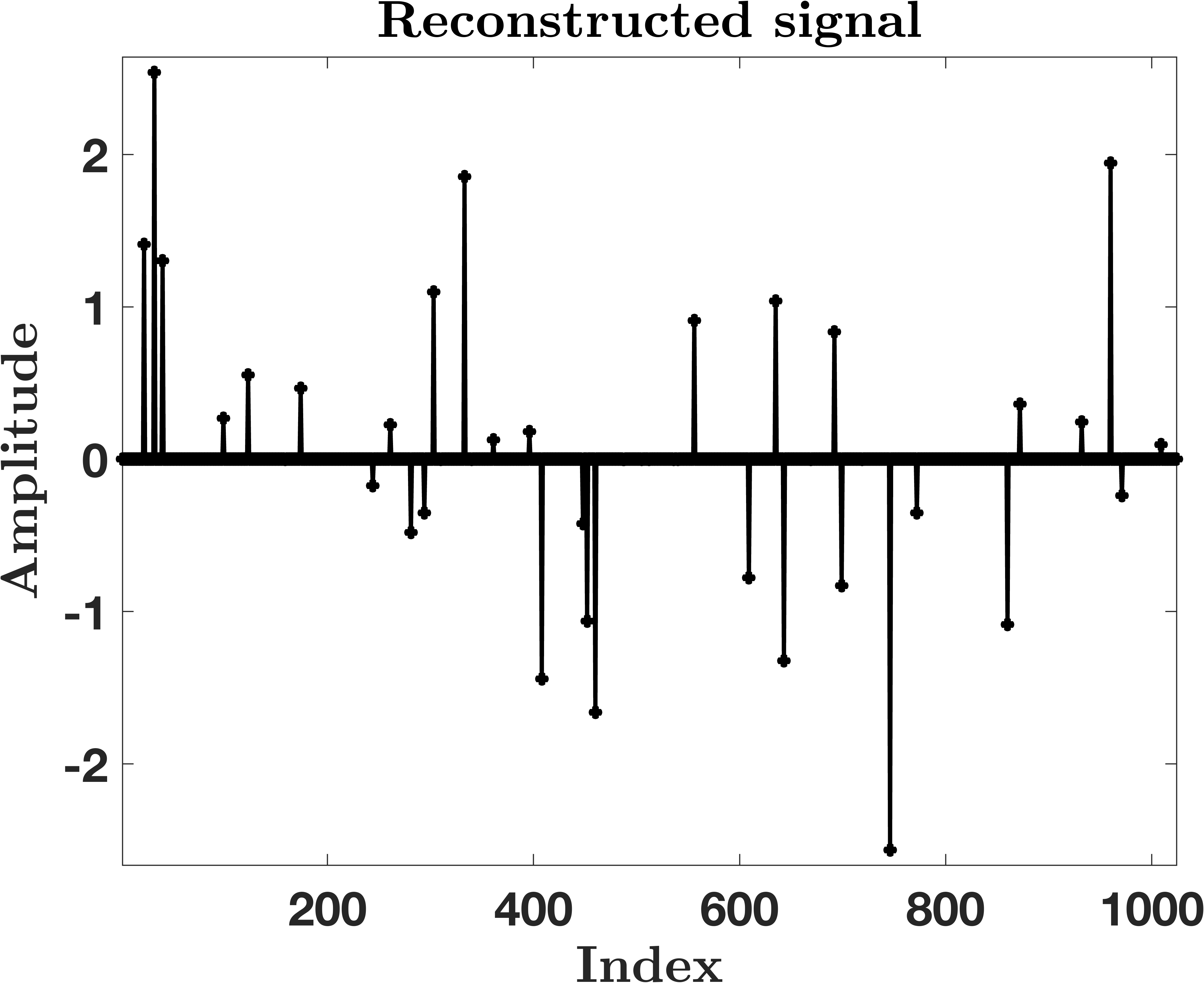}
    % \caption{Errors v.s. CPU time. }
\end{minipage}

\caption{Numerical results for Algorithm \ref{alg: restart perturbed ATPD} (ATPD-$0$) on benchmark signal reconstruction problems: (a) Sparco 9 (low resolution, size 128); (b) Sparco 10 (high resolution, size 1024); (c) Sparco 11 with $\kappa(f) = 10^8$. In each row, the panels show the error curve (left), original signal (middle), and reconstructed signal (right).}
\label{fig:ATPD-sparco-combined}
\end{figure}

We test the homotopy ATPD method (Algorithm~\ref{alg: restart perturbed ATPD}, ATPD-$0$) on benchmark signal reconstruction problems from the \textsc{Sparco} toolbox~\cite{Berg2009Sparco}. Results for Sparco 9, 10, and 11 are shown in Figure~\ref{fig:ATPD-sparco-combined}. In Algorithm~\ref{alg: restart perturbed ATPD}, we set $\theta =1$, $\varepsilon_0 = L_f$ and choose $m_0=60$. 
% We observed that the algorithm is robust with respect to $\varepsilon_0$ and $m_0$.

In the Sparco 9 and 10 problems, we reconstruct block signals generated using Heaviside matrices. Figure~\ref{fig:ATPD-sparco-combined} shows the decay of the solution error and the constraint violation, along with the original and reconstructed signals, for low-resolution recovery (size $128$) and high-resolution recovery (size $1024$). In both cases, the reconstruction error and constraint violation decrease to approximately $10^{-10}$ or below, with an approximately linear trend on the semilogarithmic plots.

For Sparco 11, the smoothed $\ell_1$ objective $f$ is less accurate at recovering the exact solution. When $\kappa(f) = 10^3$, the solution error stays around $\mathcal{O}(1)$ even when the constraint is satisfied. To improve the recovery, we increase the condition number to $\kappa(f) = 10^8$ for a better approximation of the $\ell_1$ norm. This results in the high-accuracy reconstruction shown in Figure~\ref{fig:ATPD-sparco-combined}.

\section{Conclusion}\label{sec: conclusion}

We proposed accelerated transformed primal--dual (ATPD) methods for smooth convex optimization problems with affine equality constraints. In the strongly convex regime, ATPD achieves accelerated linear convergence and matches the known lower bounds for first-order oracle complexity. In the general convex regime, the homotopy ATPD method achieves an accelerated $\mathcal O(1/N^2)$ decay of the modified Lyapunov energy after $N$ inner iterations under a uniform boundedness condition on the primal iterates at all homotopy stages.

Many problems handled by augmented Lagrangian and ADMM-type methods have nonsmooth objective functions. Extending ATPD to nonsmooth settings while preserving accelerated rates is an important direction for future work.

%\section*{Acknowledgments}
%
%The authors are sincerely grateful to the anonymous referees for their careful reading
%of the manuscript and for their valuable and constructive feedback, which has contributed significantly to improving the manuscript. The authors thank Professor Yurii Nesterov for helpful discussions at an early stage of this work, which contributed to its subsequent development and improvement.

\section*{Statements and Declarations}

\subsection*{Funding}
The first author was supported in part by National Science Foundation grants DMS-2309777 and DMS-2309785. The second author was partially supported by the School of Physical Sciences Visiting Fellows Program at the University of California, Irvine.

\subsection*{Data Availability}

The data supporting the results of this study are generated from the numerical experiments described in the manuscript. The test problems are synthetic or benchmark examples constructed using the parameters and discretization settings reported in the paper. No sensitive or personal data were used. The numerical data necessary to reproduce the reported tables and figures, together with the implementation codes, will be made available in a public repository upon publication. Additional information is available from the corresponding author upon reasonable request.

\subsection*{Competing Interests}
The authors have no competing interests to declare.

% \printbibliography
\bibliography{Optimization}  

%%%%%%%%
\end{document}